\documentclass[smallextended]{svjour3}

\usepackage[T1]{fontenc}
\usepackage{graphicx}
\usepackage{subcaption}
\usepackage{url}

\usepackage{caption}
\usepackage{amsfonts}
\usepackage{amssymb}
\usepackage{amsmath}
\usepackage{mathrsfs}
\usepackage{dsfont}
\usepackage{amsopn}

\usepackage{algorithm}
\usepackage{algorithmicx}
\usepackage{algpseudocode}

\usepackage{booktabs}
\usepackage{multirow}
\usepackage{adjustbox}
\usepackage{afterpage}

\renewcommand{\vec}[1]{\mathbf{#1}}
\newcommand{\mat}[1]{\mathbf{#1}}

\newcommand{\ie}{{\em i.e.}}

\newcommand{\inner}[1]{\langle #1 \rangle}
\newcommand{\norm}[1]{\| #1 \|}

\newcommand{\fundingname}{Funding}
\newenvironment{funding}{\par\addvspace{17pt}\small\rmfamily
\trivlist\item[\hskip\labelsep{\bfseries\fundingname}]}%
{\endtrivlist\addvspace{6pt}}

\begin{document}

\title{Accelerating a restarted Krylov method for matrix functions with randomization}
\titlerunning{Accelerating a restarted Krylov method with randomization}

\author{Nicolas L. Guidotti \and Per-Gunnar Martinsson \and Juan A. Acebr\'on \and Jos\'e Monteiro}

\institute{Nicolas L. Guidotti \at INESC-ID, Instituto Superior T\'ecnico, Universidade de Lisboa, Portugal, \\\email{nicolas.guidotti@tecnico.ulisboa.pt}  \and
Per-Gunnar Martinsson \at Department of Mathematics and Oden Institute, University of Texas at Austin, USA, \\\email{pgm@oden.utexas.edu} \and
Juan A. Acebr\'on \at Department of Mathematics, Carlos III University of Madrid, Spain, \\\email{juan.acebron@ist.utl.pt} \and
Jos\'e Monteiro \at INESC-ID, Instituto Superior T\'ecnico, Universidade de Lisboa, Portugal, \\\email{jcm@inesc-id.pt}}

\date{Received: date / Accepted: date}

\maketitle

\begin{abstract}
Many scientific applications require the evaluation of the action of a matrix function on a vector, and Krylov subspace methods are the most common ones for this task. Since the orthogonalization cost and the memory requirements can quickly become overwhelming as the basis grows, the Krylov method is often restarted after a few iterations. This paper proposes a new acceleration technique for restarted Krylov methods based on randomization. The numerical experiments show that the randomized method greatly outperforms the classical approach at the same level of accuracy. In fact, randomization can actually improve the convergence rate of restarted methods in some cases. The paper also compares the performance and stability of several existing randomized methods for solving very large, ill-conditioned problems, complementing the numerical analyses from previous studies.

\keywords{Krylov Method, Randomized algorithms, Matrix Functions, Partial Differential Equations, Finite Element Method}

\subclass{
68W20,    
65F60,    
65F50,    
65M20     
}

\begin{funding}
This work was supported by national funds through Fundação para a Ciência e a Tecnologia (FCT) under the projects URA-HPC PTDC/08838/2022 and UIDB/50021/2020 (DOI:10.54499/UIDB/50021/2020) and the grant 2022.11506.BD. NLG thanks the UT Austin Portugal program for funding his research visit to the University of Texas at Austin. JA was funded by Ministerio de Universidades and specifically the requalification program of the Spanish University System 2021-2023 at the Carlos III University. PGM recognizes support from the Office of Naval Research (N00014-18-1-2354), the National Science Foundation (DMS-2313434), and the Department of Energy ASCR (DE-SC0022251).
\end{funding}

\end{abstract}

\section{Introduction}
\label{sec:intro}

Matrix functions arise naturally in many scientific applications, for instance, in the solution of partial differential equations \cite{borner_three-dimensional_2015,guttel_rational_2013,hochbruck_exponential_2010,strutt_theory_2011}, in the analysis of complex networks \cite{benzi_matrix_2020,benzi_rational_2022,higham_functions_2008}, and in the simulation of lattice quantum chromodynamics \cite{bloch_nested_2011,van_den_eshof_numerical_2002}. Given a square matrix $\mat{A} \in \mathbb{C}^{n \times n}$, the matrix function $f(\mat{A})$ can be defined using the Cauchy integral representation \cite{higham_functions_2008}
\begin{equation}
\label{eq:cauchy_def}
f(\mathbf{A}) = \frac{1}{2\pi i} \int_{\Gamma} f(z) \; (z\mathbf{I} - \mathbf{A})^{-1} \; dz,
\end{equation}
for any function $f$ that is analytic on and inside a closed contour $\Gamma$ that encloses the spectrum $\Lambda(\mat{A})$. Some common examples include the matrix exponential $e^\mat{A}$, the matrix inverse $\mat{A}^{-1}$, and the matrix square root $\sqrt{\mat{A}}$.

Most applications are only interested in the action of $f(\mat{A})$ on a vector $\vec{b} \in \mathbb{C}^{n}$. Since explicitly forming the full matrix $f(\mat{A})$ is infeasible for large matrices, they generally employ Krylov subspace methods \cite{guttel_rational_2013,guttel_limited-memory_2020} to approximate $f(\mat{A}) \vec{b}$ directly. The key ingredient for these methods is the Arnoldi process that constructs an orthonormal basis for the Krylov subspace $\mathcal{K}_m(\mat{A}, \vec{b})$. However, the evaluation of $f(\mat{A}) \vec{b}$ requires that the entire basis be stored in memory, limiting the size of the problems that can be solved. A two-pass Lanczos method \cite{frommer_matrix_2008} can eliminate the basis storage requirements for Hermitian matrices; however, it roughly doubles the number of matrix-vector products required for constructing the basis. Furthermore, for non-Hermitian matrices, the cost of orthogonalization can quickly become overwhelming as it grows quadratically with the basis size $m$.

To overcome this challenge, restarted Krylov methods \cite{afanasjew_implementation_2008,eiermann_restarted_2006,eiermann_deflated_2011,frommer_efficient_2014,guttel_limited-memory_2020} employ a sequence of Krylov subspaces of fixed size, refining the approximation for $f(\mat{A}) \vec{b}$ at each new ``cycle''. In this manner, the program only needs to store and orthogonalize a fixed number of basis vectors at a time. However, restarted methods are often accompanied by a slower convergence rate and may even lead to stagnation or divergence \cite{schweitzer_any_2016,simoncini_convergence_2000}.

Randomization is another popular tool for reducing the orthogonalization cost in Krylov methods \cite{bergermann_deterministic_2026,burke_gmres_2025,burke_krylov_2024,cortinovis_speeding_2024,damas_randomized_2025,grigori_restoring_2026,guttel_randomized_2023,guttel_sketch-and-select_2024,jang_randomized_2024,nakatsukasa_fast_2024,palitta_sketched_2025,simoncini_stabilized_2025}. Our first contribution is to augment a restarted method with randomization. In essence, we replace the standard Arnoldi procedure with the randomized version \cite{timsit_randomized_2023} to quickly construct a non-orthogonal, but well-conditioned Krylov basis at each restart cycle. We show that this modification significantly improves the performance of the restarted methods while maintaining the same level of accuracy and stability. Randomized methods with restarts have already been proposed for exponential integrators \cite{krieger_general_2026}, linear systems \cite{burke_gmres_2025,jang_randomized_2024}, and eigendecomposition \cite{damas_randomized_2024}.

Although randomized methods perform rather well in moderately conditioned and small problems~\cite{cortinovis_speeding_2024,guttel_randomized_2023}, it is still unknown how they behave when faced with the badly conditioned matrices that are commonly found in many real-world applications. In particular, the discretization of PDEs via the finite element method (FEM) is well known to produce very large, sparse, and ill-conditioned matrices. Another example is the study of the dynamics of complex networks via the graph Laplacian. Therefore, our second contribution is an extensive study of the convergence and performance of several randomized methods through a set of large-scale numerical experiments, focusing exclusively on ill-conditioned problems.

The rest of the paper is organized as follows. Section~\ref{sec:theory} reviews the background theory of Krylov subspace methods. It also describes the randomized Arnoldi iteration and how it can be integrated into Krylov methods. Section~\ref{sec:related_works} presents other approaches for accelerating Krylov methods using randomization. Section~\ref{sec:restart_rand_krylov} describes the restarting procedure for the randomized Krylov method. Section~\ref{sec:numerical} illustrates the performance and convergence of the proposed method and compares it against the state of the art. In Section~\ref{sec:conclusion}, we conclude our paper.

\section{Theory}
\label{sec:theory}

In this section, we provide the background theory of random sketching and describe the main components of the Krylov subspace methods for evaluating $f(\mat{A}) \vec{b}$.

\subsection{Random Sketching}
\label{subsec:sketching}

For a distortion parameter $\varepsilon \in (0, 1)$, we say that the matrix $\mat{S} \in \mathbb{C}^{d \times n}$ is an $\ell_2$ embedding of a subspace $\mathbb{F} \subseteq \mathbb{C}^n$ if it satisfies
\begin{equation}
\label{eq:sketching_def1}
(1 - \varepsilon) \, \norm{\vec{x}}^2 \leq \norm{\mat{S}\vec{x} }^2 \leq (1 + \varepsilon) \, \norm{\vec{x}}^2, \quad \forall \vec{x} \in \mathbb{F},
\end{equation}
or, equivalently,
\begin{equation}
\label{eq:sketching_def2}
|\inner{\mat{S}\vec{x}, \mat{S}\vec{y}} - \inner{\vec{x}, \vec{y}}| \leq \varepsilon \, \norm{\vec{x}} \, \norm{\vec{y}}, \quad \forall \vec{x},\vec{y} \in \mathbb{F}.
\end{equation}
In practice, we do not have \textit{a priori} knowledge of $\mathbb{F}$, e.g., the Krylov subspace $\mathcal{K}_m(\mat{A}, \vec{b})$ is only available at the end of the algorithm. Therefore, we have to generate the matrix $\mat{S}$ from some random distribution that satisfies the relation (\ref{eq:sketching_def1}) with high probability. In this case, we refer to $\mat{S}$ as an \textit{oblivious subspace embedding} of $\mathbb{F}$. There are many ways to construct a subspace embedding (e.g., see \cite[Sections 8 and 9]{martinsson_randomized_2020}). Here, we focus on \textit{sparse sign matrices} \cite{cohen_nearly_2016,tropp_streaming_2019,woodruff_sketching_2014}, which take the form
\begin{equation}
\label{eq:sparse_sign_sketch}
\mat{S} = \sqrt{\frac{1}{\zeta}} \begin{bmatrix} \vec{s}_1 & \vec{s}_2 & ... & \vec{s}_n \end{bmatrix} \in \mathbb{R}^{d \times n},
\end{equation}
where $\zeta$ is a ``sparsity parameter''. The columns $\vec{s}_k \in \mathbb{R}^{d}$ are sparse random vectors with $\zeta$ nonzero entries drawn from an i.i.d.~Rademacher distribution (i.e., each entry takes $\pm 1$ with equal probability). The coordinates of the nonzero entries are uniformly chosen at random. Storing $\mat{S}$ requires $O(\zeta n)$ memory, and applying it to a vector requires $O(\zeta n)$ flops. However, it requires the use of sparse data structures and arithmetic. In terms of its quality as a random embedding, it often has similar performance to Gaussian embeddings \cite{martinsson_randomized_2020} given an appropriate choice of $\zeta$ and $d$. In this paper, we use $\zeta = 1$ (also known as \texttt{CountSketch}) and $\zeta = 4$, which were chosen to minimize the sketching cost and demonstrate the maximum attainable performance. Although we did not observe numerical issues in our experiments (Section~\ref{sec:numerical}), \texttt{CountSketch} has worse theoretical guarantees than denser sparse sign matrices and can be unreliable for certain subspaces. We recommend setting the values of $\zeta$ and $d$ based on properties of the target problem \cite{cohen_nearly_2016,li_lower_2021,martinsson_randomized_2020,tropp_streaming_2019,woodruff_sketching_2014}.

\subsection{Randomized Krylov}
\label{subsec:rand_krylov}

Recall that the Krylov subspace of order $m$ associated with $(\mat{A}, \vec{b})$ is defined as
\begin{equation}
\label{eq:krylov_def}
\mathcal{K}_m(\mat{A}, \vec{b}) = \text{span} \{\vec{b}, \mat{A} \vec{b}, \, \ldots \,, \mat{A}^{m - 1} \; \vec{b} \} \subseteq \mathbb{C}^{n}.
\end{equation}

In the classical setting, we use the Arnoldi iteration to construct an orthonormal basis $\mat{V}_m$ of $\mathcal{K}_m(\mat{A}, \vec{b})$ and an upper Hessenberg matrix $\mat{H}_m = \mat{V}_m^* \mat{A} \mat{V}_m \in \mathbb{C}^{m \times m}$ with the projection of $\mat{A}$ onto the subspace. We can then approximate $f(\mat{A})\vec{b}$ as $\vec{f}_m = \beta \mat{V}_m f(\mat{H}_m) \vec{e}_1$ with $\beta = \norm{\vec{b}}$ \cite{afanasjew_implementation_2008,eiermann_restarted_2006,eiermann_deflated_2011}.

A similar procedure can be defined for non-orthogonal bases. Let $\mat{R}_m$ be an upper Hessenberg matrix and $\mat{W}_m = [\vec{w}_1 \; \vec{w}_2 \; ... \; \vec{w}_m] \in \mathbb{C}^{n \times m}$ be a matrix whose columns determine an ascending (i.e., $\mat{W}_k = [\mat{W}_{k - 1} \; \vec{w}_k]$ with $\mat{W}_1~=~[\vec{b} / \eta]$ for a constant $\eta \in \mathbb{C}$), but not necessarily orthogonal, basis of $\mathcal{K}_m(\mat{A}, \vec{b})$. Then, we can define an \textit{Arnoldi-like decomposition} \cite{eiermann_restarted_2006} as
\begin{equation}
\label{eq:rand_arnoldi}
\mat{A} \mat{W}_m = \mat{W}_{m + 1} \mat{\underline{R}}_m =  \mat{W}_{m} \mat{R}_m + \vec{w}_{m + 1} \, r_{m + 1, m} \, \vec{e}^T_m,
\end{equation}
with
\begin{equation*}
\mat{\underline{R}}_m = \begin{bmatrix} \mat{R}_{m} \\ r_{m + 1, m} \vec{e}_m^T \end{bmatrix} \in \mathbb{C}^{(m + 1) \times m}.
\end{equation*}

\begin{lemma}[Corollary 2.2 from \cite{balabanov_randomized_2022}]
 If $\mat{S} \in \mathbb{C}^{d \times n}$ is an oblivious subspace embedding of $\mathcal{K}_m(\mat{A}, \vec{b})$, then the singular values of $\mat{W}_m$ are bounded by
\begin{equation}
\label{eq:singular_bound}
\frac{1}{\sqrt{1 + \varepsilon}} \, \sigma_{min}(\mat{S W}_m) \leq \sigma_{min}(\mat{W}_m) \leq \sigma_{max}(\mat{W}_m) \leq \frac{1}{\sqrt{1 - \varepsilon}} \, \sigma_{max}(\mat{S W}_m).
\end{equation}
\end{lemma}

Therefore, it is sufficient to orthogonalize the small sketched matrix $\mat{S W}_m$ for $\mat{W}_m$ to be well-conditioned.
This observation serves as the foundation for the Randomized Gram-Schmidt (RGS) process \cite{balabanov_randomized_2022,timsit_randomized_2023}. For each column $\vec{w}_{k + 1}$, RGS orthogonalizes the sketch $\mat{S}\vec{w}_{k + 1}$ against the sketches of the previous $k$ columns, updating the values of $\vec{w}_{k + 1}$ accordingly. This leads to a \textit{sketched-orthogonal} matrix $\mat{W}_m$, where the sketches of the columns of $\mat{W}_m$ are orthogonal among themselves. In their paper, Balabanov and Grigori \cite{balabanov_randomized_2022} analyzed the stability of the RGS process based on the distortion property of random sketching, the stability of the least-squares solver, and the finite precision of modern computers.
Algorithm~\ref{code:rand_arnoldi} describes the modified Arnoldi iteration based on the RGS process \cite{timsit_randomized_2023}.

In terms of computational complexity, forming the new basis vector $\vec{w}_{k + 1}$ requires $O(N)$ operations with $N = \mathrm{nnz}(\mat{A})$, while the cost of sketching $\vec{w}_{k + 1}$ (line 7) depends on the choice of $\mat{S}$. For a sparse sign matrix with $d = O(m)$ and $\zeta n$ nonzeros, the sketch $\mat{S} \vec{w}_{k + 1}$ can be constructed in $O(\zeta n)$ time. Orthogonalizing the sketch $\mat{S} \vec{w}_{k + 1}$ (lines~8--11) requires $O(dm) = O(m^2)$ operations and updating the vector $\vec{w}_{k + 1}$ (line 12), an additional $O(nm)$ operations. Overall, the time complexity of Algorithm~\ref{code:rand_arnoldi} is $O(Nm~+~\zeta nm~+~m^3~+~nm^2)$.

The key difference here is that at each iteration $k$, the classical Arnoldi iteration requires two passes over the basis $\mat{V}_k$, while Algorithm~\ref{code:rand_arnoldi} only needs a single pass over the full basis $\mat{W}_k$ after orthogonalizing the sketch $\mat{SW}_k$. Therefore, if $\mat{A}$ is very large and sparse, Algorithm~\ref{code:rand_arnoldi} is expected to be up to twice as fast as the standard Arnoldi procedure \cite{balabanov_randomized_2022}. Algorithm~\ref{code:rand_arnoldi} also uses a single matrix-vector product (BLAS-2) for updating $\mat{W}_k$ instead of a sequence of \texttt{dot} products and \texttt{axpy} updates (BLAS-1), which is more efficient hardware-wise.

\begin{algorithm}[t]
\caption{Randomized Arnoldi iteration for constructing the basis $\mat{W}_{m + 1}$ and the upper Hessenberg matrix $\mat{\underline{R}}_m$. Adapted from \cite[Algorithm 3]{timsit_randomized_2023}.} \label{code:rand_arnoldi}

\begin{algorithmic}[1]
\Function{RandomizedArnoldi}{$\mat{A}$, $\mat{S}$, $\vec{b}$, $m$}
\State $\alpha = \norm{\mat{S} \vec{b}}_2$
\State $\vec{u}_1 = \mat{S} \vec{b} / \alpha$ \Comment{$\mat{U} = [\vec{u}_1 ... \vec{u}_{m + 1}]$ stores $\mat{S W}$}
\State $\vec{w}_1 = \vec{b} / \alpha$
\For{$k = 1, 2,...,m$}
\State $\vec{w}_{k + 1} = \mat{A} \vec{w}_k$
\State $\vec{u}_{k + 1} = \mat{S} \vec{w}_{k + 1}$
\For{$i = 1, 2, ..., k$}
\State $r_{i, k} = \langle \vec{u}_i, \vec{u}_{k + 1}\rangle$
\State $\vec{u}_{k + 1} = \vec{u}_{k + 1} - r_{i, k} \vec{u}_i$
\EndFor
\State $\vec{w}_{k + 1} = \vec{w}_{k + 1} - \mat{W}_k \mat{R}[1:k, k]$
\State $r_{k + 1, k} = \norm{\vec{u}_{k + 1}}_2$
\State $\vec{w}_{k + 1} = \vec{w}_{k + 1} \; / \; r_{k + 1, k}$
\State $\vec{u}_{k + 1} = \vec{u}_{k + 1} \; / \; r_{k + 1, k}$
\EndFor
\State \Return $\alpha, \mat{\underline{R}}_m = [r_{i, j}], \mat{W}_{m + 1} = [\vec{w}_1 \; ... \; \vec{w}_{m + 1}]$
\EndFunction
\end{algorithmic}
\end{algorithm}

\begin{proposition}[Theorem 2.4 from \cite{eiermann_restarted_2006}]
Suppose that $\mat{W}_m$ and $\mat{R}_m$ were generated using Algorithm~\ref{code:rand_arnoldi} and that $f$ is well defined for $f(\mat{R}_m)$ and $f(\mat{A})$. Then, we can define the randomized Arnoldi approximation
$\vec{\hat{f}}_m$ for $f(\mat{A})\vec{b}$ as
\begin{equation}
\label{eq:rand_arnoldi_approx}
\vec{\hat{f}}_m = \alpha \mat{W}_m f(\mat{R}_m) \vec{e}_1.
\end{equation}
Moreover, there is a unique polynomial $q_m \in \mathcal{P}^{m - 1}$ that interpolates $f$ at the eigenvalues of $\mat{R}_m$ in the Hermite sense. Then, it holds that
\begin{equation}
\label{eq:rand_arnoldi_polynomial}
\vec{\hat{f}}_m = \alpha \mat{W}_m q_m(\mat{R}_m) \vec{e}_1 = q_m(\mat{A}) \vec{b},
\end{equation}
for any Arnoldi-like decomposition (\ref{eq:rand_arnoldi}), independently of the conditioning of $\mat{W}_m$ \cite[Corollary 2.3]{guttel_randomized_2023}.
\end{proposition}

\begin{lemma}
\label{theorem:similarity}
As both $\mat{V}_m$ and $\mat{W}_m$ are bases for the same subspace $\mathcal{K}_m(\mat{A}, \vec{b})$, there is a nonsingular, upper triangular matrix $\mat{Z} \in \mathbb{C}^{m \times m}$ such that $\alpha \mat{Z}\vec{e}_1 = \beta \vec{e}_1$ and $\mat{W}_m = \mat{V}_m \mat{Z}$. Then,
\begin{equation}
\label{eq:similarity}
\mat{Z} \mat{R}_m \mat{Z}^{-1} = \mat{H}_m - \frac{r_{m+1, m}}{z_{mm}} \, \vec{c} \, \vec{e}_m^T, \quad \vec{c} = \mat{V}_m^* \vec{w}_{m + 1}.
\end{equation}
\end{lemma}
\begin{proof}
   Replacing $\mat{W}_m$ with $\mat{V}_m \mat{Z}$ in the Arnoldi-like decomposition (\ref{eq:rand_arnoldi}), we have
   \begin{equation*}
     \mat{A} \mat{V}_m \mat{Z} = \mat{V}_m \mat{Z} \mat{R}_m + \vec{w}_{m + 1} \, r_{m + 1, m} \, \vec{e}^T_m.
   \end{equation*}
   Multiplying by $\mat{V}_m^*$ on the left and by $\mat{Z}^{-1}$ on the right,
      \begin{equation*}
     \mat{V}_m^* \mat{A} \mat{V}_m =  \mat{Z} \mat{R}_m \mat{Z}^{-1} + (\mat{V}_m^* \vec{w}_{m + 1}) \, r_{m + 1, m} \, \vec{e}^T_m \mat{Z}^{-1}.
   \end{equation*}
   By construction, $\mat{H}_m = \mat{V}_m^* \, \mat{A} \mat{V}_m$ and $\vec{e}_m^T \mat{Z}^{-1} = 1 / z_{mm} \vec{e}_m^T$ since $\mat{Z}$ is upper triangular. Setting $\vec{c} = \mat{V}_m^* \vec{w}_{m + 1}$ and rearranging the terms completes the proof.
\end{proof}

Theorem 2.2 from \cite{grigori_restoring_2026} is equivalent to Lemma~\ref{theorem:similarity} when $\vec{c} = \vec{0}$, i.e., when the last vector $\vec{w}_{m+1}$ is orthogonal to $\mat{V}_m$. With Lemma~\ref{theorem:similarity}, the randomized Arnoldi approximation $\vec{\hat{f}}_m$ can be rewritten as
\begin{equation}
\label{eq:rand_arnoldi_func_rank_one}
\vec{\hat{f}}_m = \alpha \mat{W}_m f(\mat{R}_m) \vec{e}_1 = \beta \mat{V}_m f \left(\mat{H}_m - \frac{r_{m+1, m}}{z_{mm}} \, \vec{c} \vec{e}_m^T \right) \vec{e}_1.
\end{equation}

This effectively links the approximation (\ref{eq:rand_arnoldi_func_rank_one}) with the theory in \cite{palitta_sketched_2025} that analyzes how rank-one modifications to the Hessenberg matrix $\mat{H}_m$ affect matrix functions. In particular, the numerical range $W(\mat{R}_m)$ may be significantly larger than $W(\mat{H}_m)$, and may not even be contained in $W(\mat{A})$. Yet, according to the numerical analysis from \cite{palitta_sketched_2025}, only the eigenvalues of $\mat{R}_m$ that are near or within the spectral region of $\mat{A}$ may have an influence on the convergence of the method.

\section{Related Work}
\label{sec:related_works}

A similar idea was explored in \cite{cortinovis_speeding_2024}. As the first step, their method generates a non-orthogonal basis $\mat{W}_m$ of $\mathcal{K}_m(\mat{A}, \vec{b})$ using either the incomplete Arnoldi process \cite[Chapter~6.4.2]{saad_iterative_2003} or the RGS process described in Section~\ref{subsec:rand_krylov}. Then, the projection of $\mat{A}$ onto the Krylov subspace can be computed as
\begin{equation}
\label{eq:ckn_projection}
\mat{W}_m^\dagger \mat{A} \mat{W}_m = \mat{R}_m + r_{m + 1, m} \mat{W}_{m}^\dagger \vec{w}_{m + 1} \vec{e}^T_m.
\end{equation}
For a well-conditioned basis $\mat{W}_m$ (e.g., when using Algorithm~\ref{code:rand_arnoldi}), they argue that a few iterations of LSQR \cite{paige_lsqr_1982} are sufficient to get a good enough approximation $\vec{y}_m$ for $\mat{W}_{m}^\dagger \vec{w}_{m + 1}$. If $\mat{W}_m$ is badly conditioned, they suggest using a \textit{sketch-and-precondition} approach \cite{avron_blendenpik_2010,martinsson_randomized_2020}. It is worth mentioning that solving the least-squares problem can be quite expensive for large $n$ and/or $m$. In both cases, $f(\mat{A}) \vec{b}$ can be approximated as
\begin{equation}
\label{eq:ckn_matfunc}
\vec{f}_{m}^{ckn} = \eta \mat{W}_m f(\mat{Y}_m) \, \vec{e}_1,
\end{equation}
with $\mat{Y}_m = \mat{W}_m^\dagger \mat{A} \mat{W}_m = \mat{R}_m + r_{m + 1, m} \, \vec{y}_m \, \vec{e}^T_m$ and $\eta = \norm{\vec{b}}$ for the incomplete orthogonalization or $\eta = \norm{\mat{S} \vec{b}}$ for the randomized Arnoldi (Algorithm~\ref{code:rand_arnoldi}). The expression (\ref{eq:rand_arnoldi_approx}) is a simplification of (\ref{eq:ckn_matfunc}) since generally $\mat{W}_m^\dagger \mat{A} \mat{W}_m \neq \mat{R}_m$ (see Lemma~\ref{theorem:similarity}). In practice, the upper Hessenberg matrices $\mat{Y}_m$ and $\mat{R}_m$ only differ in the last column, such that this difference may have only a minor influence on the first column of $f(\mat{W}_m^\dagger \mat{A} \mat{W}_m)$. A similar approximation was proposed in \cite[Section 3.2]{cortinovis_speeding_2024}, but for the incomplete orthogonalization case. Likewise, restarted methods often use a similar approximation with a non-orthonormal basis $\mat{W}_m$; see \cite{eiermann_restarted_2006} for more details.

The authors of \cite{grigori_restoring_2026} modify Algorithm~\ref{code:rand_arnoldi} to produce a Hessenberg matrix $\mat{R_m}$ that is similar to the one generated by the classical Arnoldi iteration. This is achieved by first constructing a sketch-orthogonal basis $\mat{W}_m$ and then solving a least-squares problem to orthogonalize the last basis vector against all previous vectors.

In \cite{guttel_randomized_2023}, the authors propose to first construct a non-orthogonal Krylov basis $\mat{W}_m$ using an incomplete Arnoldi process. They then approximately orthogonalize $\mat{W}_m$, working only with the sketch of the basis via a thin QR decomposition (``\textit{basis whitening}'' \cite{guttel_randomized_2023,nakatsukasa_fast_2024}). They call this method \texttt{sFOM}, which was later refined in \cite{palitta_sketched_2025}. It is worth mentioning that in the same paper \cite{guttel_randomized_2023}, the authors also propose another randomized Krylov method --- \texttt{sGMRES} --- which requires numerical quadrature rules for evaluating matrix functions.

For solving ODEs via exponential integrators, the authors of \cite{krieger_general_2026} derived an \textit{a posteriori} error estimation formula for \texttt{sFOM} based on the residual of the underlying differential equation. They then use this information to monitor the convergence and decide when to restart. The restart is done by adaptively changing the length of the time step taken during the integration (residual-time restart). The authors include a direct comparison against the restarted randomized method from this paper.

Finally, a recent survey \cite{damas_randomized_2025} summarizes the different randomization techniques that are applicable to Krylov methods, focusing primarily on the basis construction methods.

\section{Restarted Randomized Krylov}
\label{sec:restart_rand_krylov}

Instead of using a single, large Krylov subspace for computing $f(\mat{A}) \vec{b}$, restarted methods \cite{afanasjew_implementation_2008,eiermann_restarted_2006,guttel_limited-memory_2020} employ a sequence of Krylov subspaces of fixed size, refining the approximation $\vec{f}_m$ at each new ``cycle''.

In this section, we propose a restarting procedure for the randomized Krylov method from Section~\ref{subsec:rand_krylov}. Let us consider the first two restart cycles. Suppose that after $m$ iterations of Algorithm~\ref{code:rand_arnoldi}, we obtain the following decomposition
\begin{equation*}
\mat{A} \mat{W}_m^{(1)} = \mat{W}_{m}^{(1)} \mat{R}_m^{(1)} + \vec{w}_{m + 1}^{(1)} \, r_{m + 1, m}^{(1)} \, \vec{e}^T_m,
\end{equation*}
with $\vec{w}_1^{(1)} = \vec{b} / \alpha$ and the approximation $\vec{\hat{f}}_1 = \alpha \mat{W}_m^{(1)} \; f(\mat{R}_m^{(1)}) \vec{e}_1$. We then ``restart'' the randomized Arnoldi process, now with $\vec{w}_{m + 1}^{(1)}$ as the starting vector, and obtain a second decomposition
\begin{equation*}
\mat{A} \mat{W}_m^{(2)} = \mat{W}_{m}^{(2)} \mat{R}_m^{(2)} + \vec{w}_{m + 1}^{(2)} \, r_{m + 1, m}^{(2)} \, \vec{e}^T_m.
\end{equation*}
Combining both decompositions, we have
\begin{equation}
\label{eq:restart_combined}
\mat{A} \mat{W}_{2m} =  \mat{W}_{2m} \mat{R}_{2m} + \vec{w}_{m + 1}^{(2)} \, r_{m + 1, m}^{(2)} \, \vec{e}^T_{2m},
\end{equation}
where
\begin{equation*}
\mat{R}_{2m} = \begin{bmatrix} \mat{R}_{m}^{(1)} & \mat{0} \\ r_{m + 1, m}^{(1)} \vec{e}_1\vec{e}_m^T & \mat{R}_{m}^{(2)}  \end{bmatrix} \quad \text{and} \quad \mat{W}_{2m} = [\mat{W}_m^{(1)} \quad \mat{W}_m^{(2)}].
\end{equation*}
The Krylov approximation $\vec{\hat{f}}_2$ associated with (\ref{eq:restart_combined}) is then defined as
\begin{equation}
\label{eq:restart_approx}
\vec{\hat{f}}_{2} = \alpha \; \mat{W}_{2m} \; f(\mat{R}_{2m}) \; \vec{e}_1.
\end{equation}
Due to the block triangular structure of $\mat{R}_{2m}$, $f(\mat{R}_{2m})$ has the following form
\begin{equation*}
f(\mat{R}_{2m}) = \begin{bmatrix} f(\mat{R}_{m}^{(1)}) & \mat{0} \\ \mat{F}_{2, 1}  & f(\mat{R}_{m}^{(2)})  \end{bmatrix},
\end{equation*}
and thus, (\ref{eq:restart_approx}) can be rewritten as
\begin{equation}
\label{eq:restart_approx_update}
\vec{\hat{f}}_{2} = \vec{\hat{f}}_1 + \alpha \; \mat{W}_{m}^{(2)} \; \mat{F}_{2,1} \; \vec{e}_1.
\end{equation}
Therefore, as long as we can compute $\mat{F}_{2,1}$, we can update the Krylov approximation $\vec{\hat{f}}_{2}$ without needing to store the basis from the previous cycle. For the method to be numerically stable \cite{eiermann_restarted_2006,guttel_limited-memory_2020}, the matrix $\mat{F}_{2,1}$ is taken directly from the bottom-left block of $f(\mat{R}_{2m})$, which in turn entails the computation of the full matrix $f(\mat{R}_{2m})$. Algorithm~\ref{code:restart_krylov} describes the complete restarting procedure for an arbitrary number of cycles. Note that the size of $\mat{R}_{km}$ grows with each restart, such that evaluating $f(\mat{R}_{km})$ may become too expensive to compute after a large number of restarts \cite{eiermann_restarted_2006,guttel_limited-memory_2020}. In this scenario, additional acceleration techniques may be required, such as spectral deflation \cite{eiermann_deflated_2011}, thick restarting \cite{baker_technique_2005}, and weighted inner products \cite{embree_weighted_2017}. Note that the Krylov method proposed by \cite{frommer_efficient_2014,guttel_limited-memory_2020} does not have this limitation, but they rely on quadrature rules for evaluating Cauchy-Stieltjes functions and the exponential; hence, they are outside the scope of this paper.

\begin{algorithm}[t]
\caption{Restarted Krylov method for evaluating $f(\mat{A})\vec{b}$ based on Algorithm~\ref{code:rand_arnoldi}. $tol$ is the tolerance and $k_{max}$ is the maximum number of cycles.} \label{code:restart_krylov}

\begin{algorithmic}[1]
\Function{RandomizedRestartedKrylov}{$\mat{A}$, $\mat{S}$, $\vec{b}$, $m$, $tol$, $k_{max}$}
\State $\alpha, \underline{\mat{R}}_m^{(1)}, \mat{W}_{m + 1}^{(1)} = $ \Call{RandomizedArnoldi}{$\mat{A}$, $\mat{S}$, $\vec{b}$, $m$}
\State $\vec{y} = \alpha \mat{W}_m^{(1)} \; f(\mat{R}_m^{(1)}) \vec{e}_1$
\State $\vec{\hat{f}}_1 = \vec{y}$
\For{$k = 2, ..., k_{max}$}
\State \textbf{if} $\norm{\vec{y}} \leq tol$ \textbf{then break}
\State $\sim, \underline{\mat{R}}_m^{(k)}, \mat{W}_{m + 1}^{(k)} = $ \Call{RandomizedArnoldi}{$\mat{A}$, $\mat{S}$, $\vec{w}_{m + 1}^{(k - 1)}$, $m$}
\State $\mat{R}_{km} = \begin{bmatrix} \mat{R}_{(k - 1) m} & \mat{0} \\ r_{m + 1, m}^{(k - 1)} \vec{e}_1\vec{e}_{(k - 1)m}^T & \mat{R}_{m}^{(k)}  \end{bmatrix}$
\State $\mat{F} = f(\mat{R}_{km})$
\State $\vec{y} = \alpha \; \mat{W}_{m}^{(k)} \; \mat{F}[(k - 1) m + 1: km, \, 1]$
\State $\vec{\hat{f}}_k = \vec{\hat{f}}_{k - 1} + \vec{y}$
\EndFor
\State \Return $\vec{\hat{f}}_k$
\EndFunction
\end{algorithmic}
\end{algorithm}

As with other iterative procedures, we need to estimate the error $\norm{f(\mat{A})\vec{b} - \vec{\hat{f}}_k}$ in order to determine when the restarted algorithm should stop. A simple error estimate that is often used in Krylov methods is the difference between two successive restart cycles:
\begin{equation*}
\norm{f(\mat{A})\vec{b} - \vec{\hat{f}}_k} \approx \norm{\vec{\hat{f}}_{k - 1} - \vec{\hat{f}}_k}.
\end{equation*}

The convergence of (classical) restarted methods has been established for entire functions of order $1$ \cite[Theorem 4.2]{eiermann_restarted_2006} and Stieltjes functions \cite{frommer_convergence_2014}, while the general case remains an open problem. The convergence of restarted methods in the context of Laplace transforms for Hermitian matrices is discussed in \cite{tsolakis_efficient_2023} (our focus is on the non-symmetric case). With randomization, a rigorous analysis becomes even more challenging. Although \cite{cortinovis_speeding_2024,guttel_randomized_2023,palitta_sketched_2025} provide some theoretical results, the convergence of unrestarted randomized methods is still not fully understood due to many obstacles in the analysis, such as the fact that the numerical range of the Hessenberg matrix $\mat{R}_{m}$ is difficult to control and may not be contained in the numerical range of $\mat{A}$ \cite{grigori_restoring_2026,palitta_sketched_2025}. Considering these limitations, we will analyze the convergence of the proposed algorithm in terms of numerical experiments.

\section{Numerical Experiments}
\label{sec:numerical}

To evaluate the performance and stability of the randomized Krylov methods, we solve a set of linear partial differential equations (PDEs) that arise in different scientific applications. Following the method of lines~\cite{schiesser_numerical_1991}, we first discretize the spatial variables of the PDE, transforming the original problem into a system of coupled ordinary differential equations with time as the independent variable. The initial value problem can then be solved by evaluating some function $f$ over the coefficient matrix~$\mat{A}$. Another example consists of simulating the diffusion on a network using the graph Laplacian. In summary, we analyze the following methods in this section:

\begin{itemize}
  \item \texttt{arnoldi}: Classical method without restarts.
  \item \texttt{incomplete}: Classical method with an incomplete orthogonalization \cite{saad_iterative_2003}.
  \item \texttt{rand}: Randomized method from Section~\ref{subsec:rand_krylov}.
  \item \texttt{rand-ls}: Randomized method proposed in \cite[Algorithm 3.1]{cortinovis_speeding_2024}. The Krylov basis is generated using Algorithm~\ref{code:rand_arnoldi}. The least-squares problem is solved with LSMR \cite{fong_lsmr_2011} with a tolerance of $\epsilon = 10^{-8}$ and a maximum of $20$ iterations.
  \item \texttt{sFOM}: Randomized method proposed in \cite{guttel_randomized_2023,palitta_sketched_2025}.
  \item \texttt{restart}: Classical restarted method from \cite{afanasjew_implementation_2008,eiermann_restarted_2006,guttel_limited-memory_2020}.
  \item \texttt{deflated}: Classical restarted method with spectral deflation from \cite{eiermann_deflated_2011}. In each restart cycle, we deflate the $\ell=5$ smallest Ritz values from the previous cycle.
  \item \texttt{restart-rand}: Randomized restarted method from Section~\ref{sec:restart_rand_krylov}.
\end{itemize}

We use a $d \times n$ sparse sign matrix with the sparsity parameter $\zeta$ as the sketching matrix $\mat{S}$ in all randomized methods. All algorithms were implemented in \texttt{C++} using Intel MKL v2026.0 for all BLAS/LAPACK routines. The code was compiled with \texttt{clang} v22.1.8. The results shown here are the arithmetic mean of five runs with different seeds, and the $\pm$ values reported in the tables are the standard deviation. The range of values for each curve is indicated by the shaded area. Note that the initial conditions are always the same across all runs. The observed distortion $\varepsilon_{obs}$ was calculated as $\varepsilon_{obs} = \max\left( \sigma_{min}^{-2}(\mat{W}_m) - 1,\;\; 1 - \sigma_{max}^{-2}(\mat{W}_m) \right)$.

\subsection{Convection-Diffusion}
\label{subsec:conv_diff}

The first example consists of solving a three-di\-men\-sional convection-diffusion equation over a domain $\Omega$
\begin{equation}
\label{eq:conv_diff_def}
\frac{\partial}{\partial t} u(\vec{x}, t) = \nabla(\alpha \nabla u(\vec{x}, t)) + \beta \nabla u(\vec{x}, t) + g(\vec{x}, t),
\end{equation}
for a time $t > 0$, the space variable $\vec{x} = (x, y, z) \in \Omega$, a diffusion coefficient $\alpha$, and a convection coefficient $\beta$. In this example, we consider $\alpha$ and $\beta$ to be constant in both time and space. For the problem to be well defined, we set the initial condition as $u(\vec{x}, 0) = u_0(\vec{x})$ and impose boundary conditions on $\partial \Omega$. After applying the standard Galerkin finite element procedure \cite{larson_finite_2013,zienkiewicz_finite_2013} over the domain $\Omega$, we obtain the following system of equations:
\begin{equation}
\label{eq:conv_diff_fem}
\mat{M} \frac{d}{dt} \vec{\hat{u}} = -\mat{A} \vec{\hat{u}} + \mat{M} \vec{g}(t),
\end{equation}
where $\mat{M}$ is the mass matrix, $\mat{A} = \alpha \mat{D} + \beta \mat{C}$ is the stiffness matrix, and $\vec{g}(t)$ is the load vector. Here, $\mat{D}$ and $\mat{C}$ correspond to the matrices related to the discretization of the diffusion and convection parts in the equation, respectively. We can then write the solution of (\ref{eq:conv_diff_fem}) as
\begin{equation}
\label{eq:conv_diff_sol}
\vec{\hat{u}}(t) = \exp(-\mat{L} \, t) \, \vec{\hat{u}}_0 + \int_{0}^{t} \exp(-\mat{L} \, s) \, \vec{g}(t - s) \, ds,
\end{equation}
with $\mat{L} = \mat{M}^{-1}\mat{A}$. To simplify the computation, we assume that the mass matrix $\mat{M}$ is lumped \cite{larson_finite_2013,zienkiewicz_finite_2013}, such that all of its mass is concentrated on the diagonal. We also assume that the load vector $\vec{g}$ remains constant throughout the entire simulation, such that the solution (\ref{eq:conv_diff_sol}) can be simplified to
\begin{equation}
\label{eq:conv_diff_sol2}
\vec{\hat{u}}(t) = \exp(-\mat{L} \, t) \, \vec{\hat{u}}_0 - \mat{L}^{-1} \exp(-\mat{L} \, t) \vec{g} + \mat{L}^{-1} \vec{g},
\end{equation}
or, equivalently,
\begin{equation}
\label{eq:conv_diff_sol3}
\vec{\hat{u}}(t) = \vec{\hat{u}}_0 + t\varphi_1(-\mat{L} t) \, \vec{b}, \quad \varphi_1(z) = \frac{e^z - 1}{z},
\end{equation}
with $\vec{b} = \vec{g} - \mat{L} \vec{\hat{u}}_0$. The function $\varphi_1(z)$ is known as the ``phi function'' in the exponential integrator literature \cite{hochbruck_exponential_2010}. It is worth mentioning that evaluating (\ref{eq:conv_diff_sol3}) is more stable and faster than (\ref{eq:conv_diff_sol2}) as it avoids solving a linear system with $\mat{L}$, which can be quite problematic when $\mat{L}$ has an eigenvalue near the origin. This is not an issue when using (\ref{eq:conv_diff_sol3}) as $\varphi_1(z)$ is an entire function.

\begin{figure}[t]
\centering
    \begin{subfigure}{0.3 \textwidth}
        \centering
        \includegraphics[width=\linewidth]{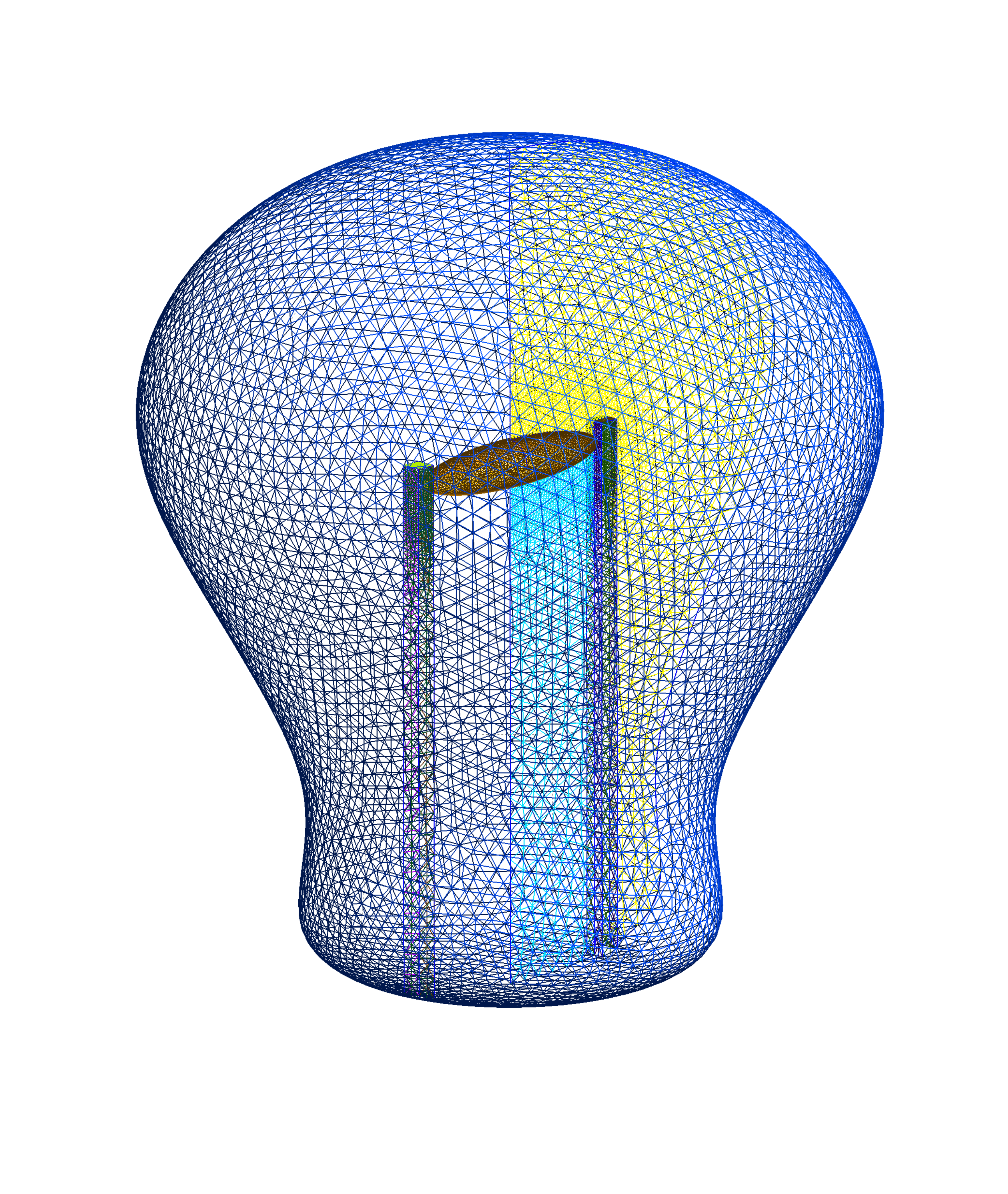}
    \end{subfigure}
    \begin{subfigure}{0.3  \textwidth}
        \centering
        \includegraphics[width=\linewidth]{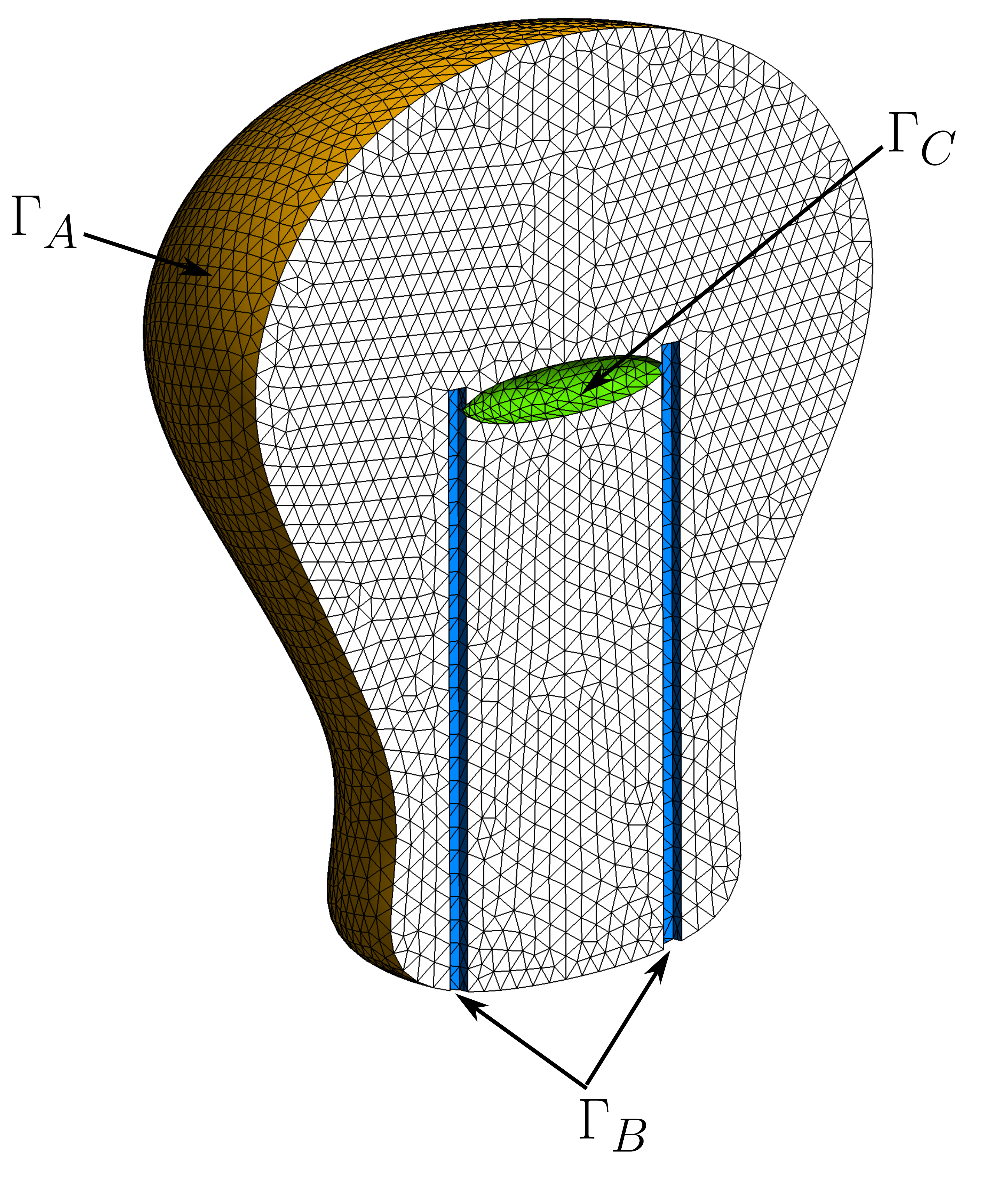}
    \end{subfigure}
    \caption{The finite element mesh and its boundary conditions for the convection-diffusion example.}
    \label{fig:fem_geo}
\end{figure}

Fig.~\ref{fig:fem_geo} shows the geometry of the domain used in this numerical experiment. The object has a finer discrete mesh near the boundary $\Gamma_C$, with the finite element size increasing towards the boundary $\Gamma_A$. The discrete mesh was generated with \texttt{gmsh} \cite{geuzaine_gmsh_2009}, while the mass and stiffness matrices were assembled using \texttt{FreeFem++}~\cite{hecht_new_2012} and P1 finite elements. We consider the following boundary conditions:
\begin{subequations}
    \begin{align}
      & u(\vec{x}, t) = 0 \quad \text{at} \quad x \in \Gamma_A, \label{eq:lamp_boundary_a} \\
      & \nabla u(\vec{x}, t) \cdot \vec{n} = 0 \quad \text{at} \quad x \in \Gamma_B,  \label{eq:lamp_boundary_b} \\
      & u(\vec{x}, t) = 1 \quad \text{at} \quad x \in \Gamma_C. \label{eq:lamp_boundary_c}
    \end{align}
\end{subequations}
Here, $\vec{n}$ denotes the outward normal vector to the boundary $\Gamma_B$. Let us assume that the rows of $\mat{L}$ are ordered such that the first $n_{i}$ rows correspond to the nodes in the interior of $\Omega$ or on the boundary $\Gamma_B$, while the remaining $n_{b}$ rows correspond to the nodes at the boundaries $\Gamma_A$ and $\Gamma_C$. Then, to impose the boundary conditions (\ref{eq:lamp_boundary_a}) and (\ref{eq:lamp_boundary_c}), we set the matrix $\mat{L}$, the load vector $\vec{g}$, and the initial conditions $\vec{u}_0$ as
\begin{equation*}
\mat{L} =
\begin{bmatrix}
\mat{L}_{11} & \mat{L}_{12} \\
\mat{0} & \mat{I}
\end{bmatrix}
\qquad
\vec{g}(\vec{x}) =
\begin{cases}
1, \quad \vec{x} \in \Gamma_C \\
0, \quad \text{otherwise}
\end{cases}
\qquad
\vec{u}_0 (\vec{x}) =
\begin{cases}
0, \quad \vec{x} \in \Gamma_A \\
1, \quad \vec{x} \in \Gamma_C \\
\xi, \quad \text{otherwise}
\end{cases}
\end{equation*}
where $\mat{L}_{11}$ is the upper-left $n_{i} \times n_{i}$ block from $\mat{L}$, $\mat{L}_{12}$ is the upper-right $n_{i} \times n_{b}$ block from $\mat{L}$, $\mat{I}$ is the $n_{b} \times n_{b}$ identity matrix, and $\xi$ is a random number drawn from the Gaussian distribution $\mathcal{N}(0.5, 0.25)$. The Neumann boundary condition (\ref{eq:lamp_boundary_b}) was enforced during the assembly of the matrix $\mat{A}$.

\begin{table}[t]
\centering
\caption{Extremal eigenvalues and condition number of $\mat{L}$ for a mesh with a finite element size between $0.01$ and $0.05$.}
\label{tab:conv-diff-eigen}
\begin{tabular}{@{}llllrlrr@{}}
\toprule
         & $\alpha$ & $\beta$ & $\lambda_{min}$ & $\lambda_{max}$ & $\sigma_{min}$ & $\sigma_{max}$ & $\kappa_2(\mat{L})$ \\ \midrule
Test (a) & $0.1$    & $0.01$  & $0.120$         & $23,363$        & $0.0706$         & $25,842$        & $365,898$           \\ \midrule
Test (b) & $0.01$   & $0.01$  & $0.0194$        & $2,337$         & $0.00813$         & $2,585$ & $317,921$           \\ \midrule
Test (c) & $0.01$   & $1$     & $1$             & $2,437$         & $0.0283$         & $2,845$ & $100,644$           \\ \bottomrule
\end{tabular}
\end{table}

\begin{figure}[t]
  \centering
  \includegraphics[width=\linewidth]{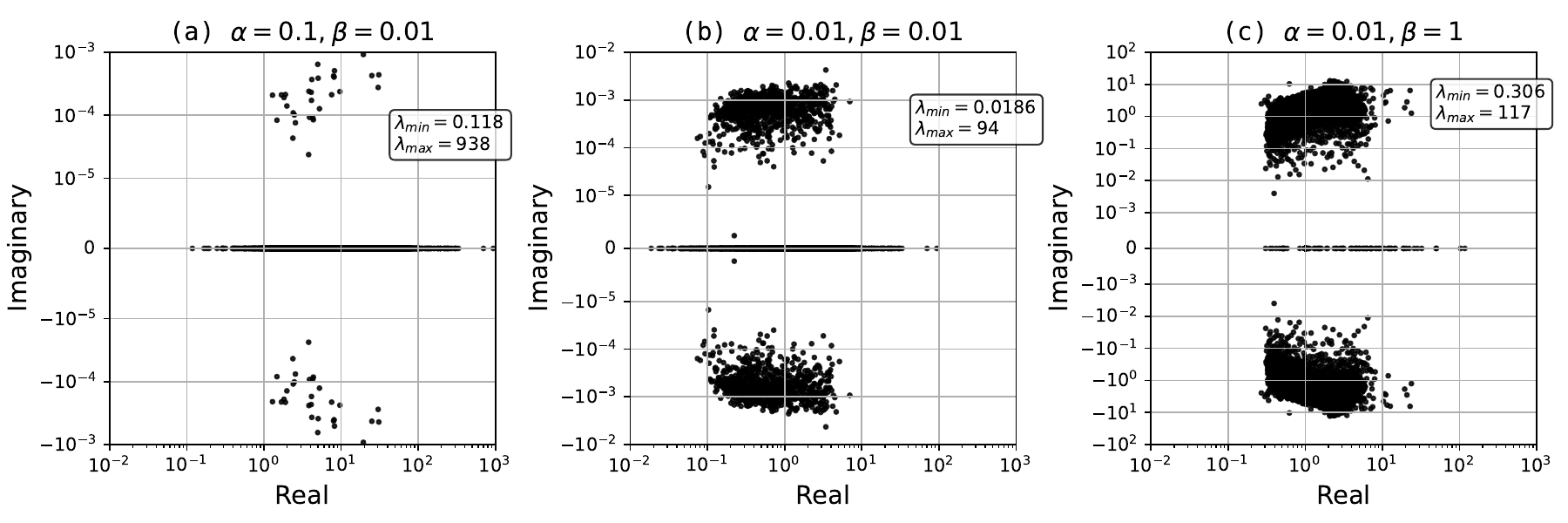}
  \caption{The spectrum of the matrix $\mat{L}$ for a coarse mesh with the size of the finite elements ranging between $0.12$ and $0.6$.}
  \label{fig:fem_eigen}
\end{figure}

To gain some insight into the spectral properties of $\mat{L}$, we set the discrete mesh to be very coarse, such that the resulting matrix $\mat{L}$ is sufficiently small for its full eigendecomposition to be feasible. As shown in Fig.~\ref{fig:fem_eigen}, the entire spectrum of $\mat{L}$ is located in the right half-plane, which indicates that the solution (\ref{eq:conv_diff_sol3}) will converge to a steady state after a sufficiently long time $t$ has passed. If $\alpha \geq \beta$ (i.e., the diffusion term is equal to or greater than the convection term), the majority of the eigenvalues of $\mat{L}$ lie on or near the real axis. If the dynamics of the problem are dominated by the convection (i.e., $\beta \gg \alpha$), most eigenvalues are complex with a large imaginary part. Generally speaking, the spectrum is wider for higher values of $\alpha$.

For the remaining numerical examples in this section, we consider a finite element size between $0.01$ and $0.05$, generating a matrix $\mat{L}$ with $4,801,565$ rows and $73,074,426$ nonzeros. Although the matrix is too large to compute the full eigendecomposition, we can estimate the extremal eigenvalues and singular values of $\mat{L}$ using ARPACK~\cite{lehoucq_arpack_1998}. They are displayed in Table~\ref{tab:conv-diff-eigen}. We adopted as reference the solution obtained with the \texttt{restart} method considering a tolerance of $3 \times 10^{-12}$ and a restart length $m_r = 100$. The matrix phi-function $\varphi_1$ was calculated using the scaling-and-squaring method from~\cite{al-mohy_new_2010,higham_scaling_2005}.

\begin{table}[t]
\centering
\caption{The relative $\ell_2$ error, condition number of the basis $\kappa_2(\mat{W}_m)$, and observed distortion $\varepsilon_{obs}$ at the last iteration of non-restarted randomized methods. The sketching dimension $d$ was set to $3200$. For \texttt{sFOM}, we report $\kappa_2(\mat{W}_m \mat{T}_m^{-1})$ instead of $\kappa_2(\mat{W}_m)$.}
\begin{subtable}{\textwidth}
\centering
\caption{$\alpha = 0.1, \beta = 0.01, t = 1, m = 1000$}
\label{tab:conv_diff_test_a_full}
\begin{tabular}{ccccc}
\hline
method & $\zeta$ & $\frac{\norm{\hat{\vec{u}} - \vec{u}_{ref}}}{\norm{\vec{u}_{ref}}}$ & $\kappa_2(\mat{W}_m)$ & $\varepsilon_{obs}$ \\
\hline
\texttt{rand} & $1$ & $(4.28 \pm 3.87) \times 10^{-13}$ & $8.32 \pm 3.18$ & $1.80 \pm 0.22$ \\
\texttt{rand} & $4$ & $(3.87 \pm 3.33) \times 10^{-13}$ & $3.51 \pm 0.01$ & $1.43 \pm 0.02$ \\
\texttt{sFOM} & $4$ & $(3.98 \pm 2.69) \times 10^{-13}$ & $70.31 \pm 9.03$ & $226.66 \pm 31.46$ \\
\texttt{rand-ls} & $4$ & $(3.99 \pm 3.61) \times 10^{-13}$ & $3.51 \pm 0.01$ & $1.43 \pm 0.02$ \\
\hline
\end{tabular}
\end{subtable} \\
\vspace{1.5em}
\begin{subtable}{\textwidth}
\caption{$\alpha = 0.01, \beta = 0.01, t = 1, m = 300$}
\label{tab:conv_diff_test_b_full}
\centering
\begin{tabular}{ccccc}
\hline
method & $\zeta$ & $\frac{\norm{\hat{\vec{u}} - \vec{u}_{ref}}}{\norm{\vec{u}_{ref}}}$ & $\kappa_2(\mat{W}_m)$ & $\varepsilon_{obs}$ \\
\hline
\texttt{rand} & $1$ & $(3.06 \pm 2.05) \times 10^{-14}$ & $2.06 \pm 0.26$ & $0.75 \pm 0.08$ \\
\texttt{rand} & $4$ & $(3.63 \pm 2.32) \times 10^{-14}$ & $1.89 \pm 0.01$ & $0.72 \pm 0.02$ \\
\texttt{sFOM} & $4$ & $(6.74 \pm 5.97) \times 10^{-14}$ & $32.93 \pm 2.15$ & $412.31 \pm 39.05$ \\
\texttt{rand-ls} & $4$ & $(3.91 \pm 2.35) \times 10^{-14}$ & $1.89 \pm 0.01$ & $0.72 \pm 0.02$ \\
\hline
\end{tabular}
\end{subtable} \\
\vspace{1.5em}
\begin{subtable}{\textwidth}
\centering
\caption{$\alpha = 0.01, \beta = 1, t = 1, m = 500$}
\label{tab:conv_diff_test_c_full}
\begin{tabular}{ccccc}
\hline
method & $\zeta$ & $\frac{\norm{\hat{\vec{u}} - \vec{u}_{ref}}}{\norm{\vec{u}_{ref}}}$ & $\kappa_2(\mat{W}_m)$ & $\varepsilon_{obs}$ \\
\hline
\texttt{rand} & $1$ & $(2.34 \pm 1.90) \times 10^{-14}$ & $3.44 \pm 0.84$ & $1.07 \pm 0.10$ \\
\texttt{rand} & $4$ & $(2.95 \pm 2.23) \times 10^{-14}$ & $2.30 \pm 0.01$ & $0.94 \pm 0.01$ \\
\texttt{sFOM} & $4$ & $(7.36 \pm 6.78) \times 10^{-14}$ & $30.86 \pm 3.35$ & $148.68 \pm 7.15$ \\
\texttt{rand-ls} & $4$ & $(2.59 \pm 1.25) \times 10^{-14}$ & $2.30 \pm 0.01$ & $0.94 \pm 0.01$ \\
\hline
\end{tabular}
\end{subtable}
\label{tab:conv_diff_full}
\end{table}

\begin{figure}[t]
  \centering
  \includegraphics[width=\textwidth]{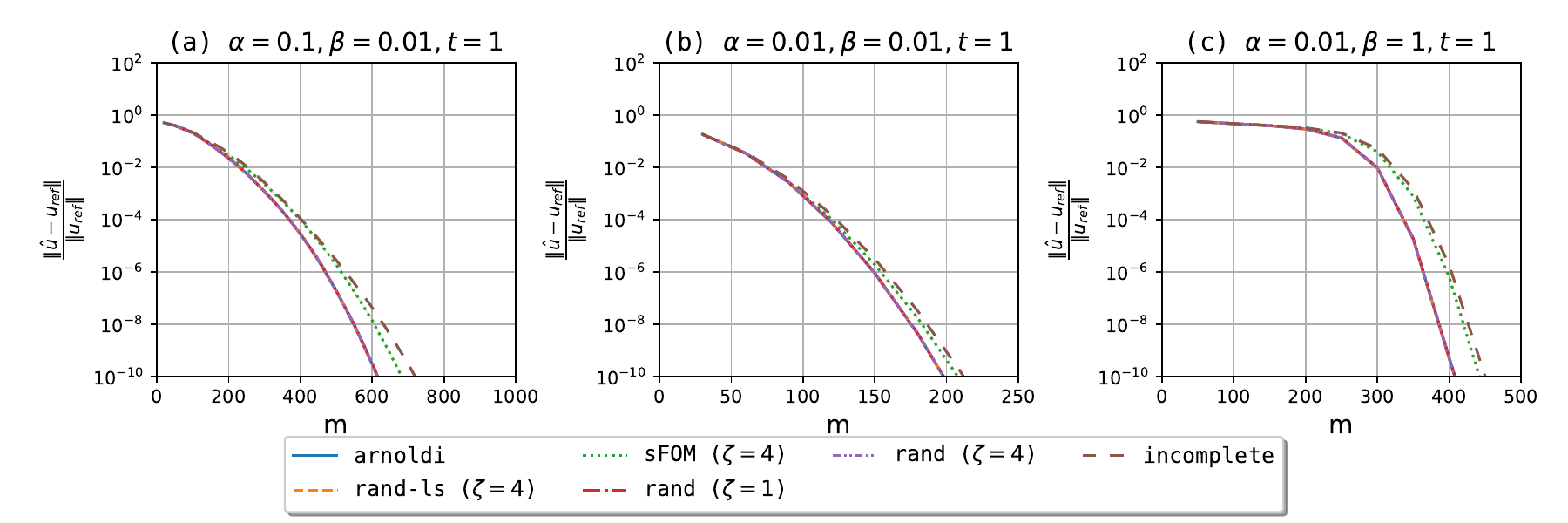}
  \caption{Convergence curves for non-restarted methods for the convection-diffusion problem. The truncation parameter $k$ was set to $20$ and the sketching dimension $d$ was set to $3200$.}
  \label{fig:conv_diff_convergence}
\end{figure}

\begin{figure}[t]
	\centering	\includegraphics[width=\textwidth]{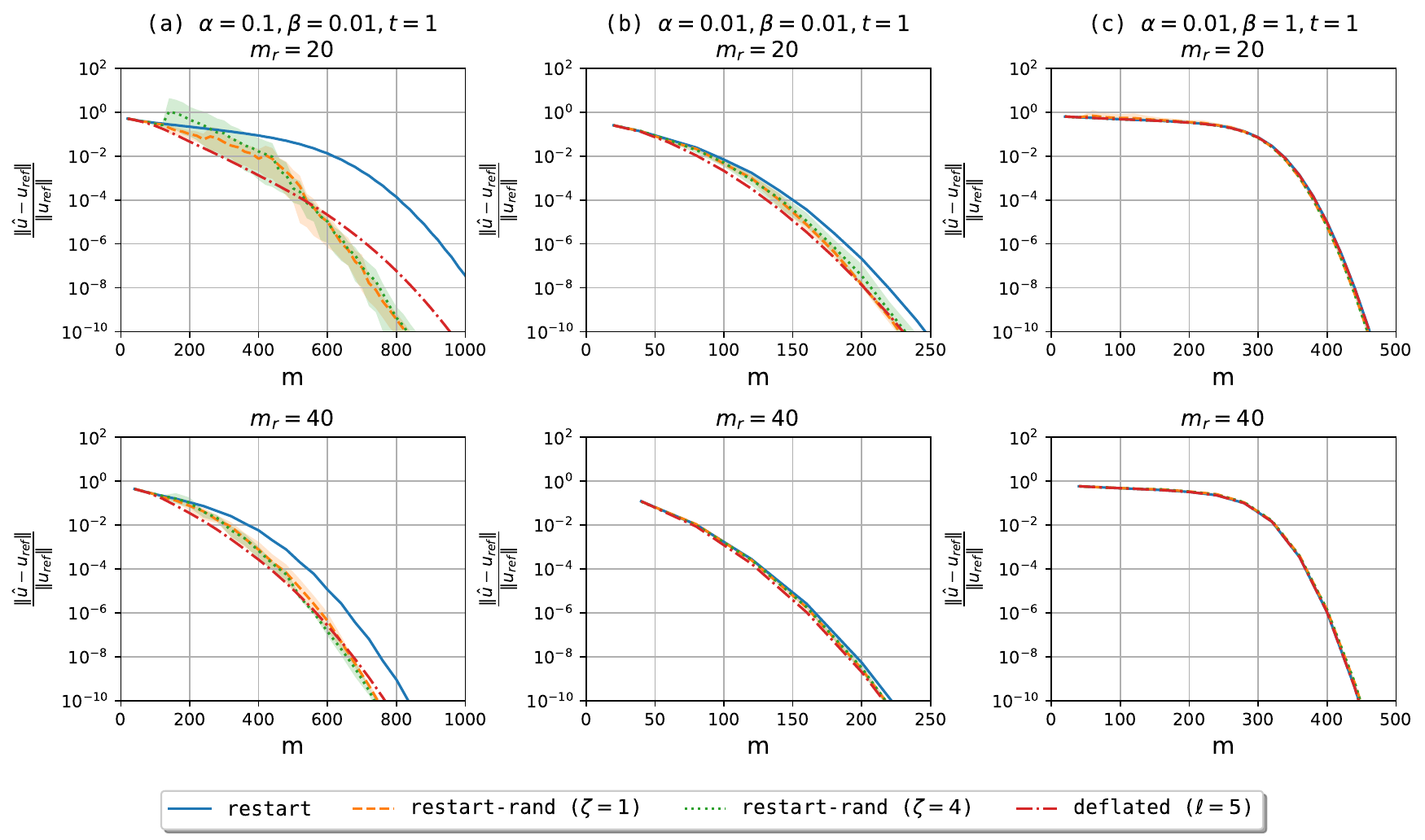}
	\caption{Convergence curves for restarted methods for the convection-diffusion example. The sketching dimension $d$ was set to $320$.}
	\label{fig:conv_diff_restart}
\end{figure}

Fig.~\ref{fig:conv_diff_convergence} shows the convergence curves for non-restarted methods, while Table~\ref{tab:conv_diff_full} displays some key information for each randomized method. From Table~\ref{tab:conv-diff-eigen}, we can infer that the spectrum of $\mat{L}$ in test~(a) is significantly wider than in the other tests, and thus it requires a much larger $m$ to approximate the entire spectrum. Likewise, assuming that the spectrum is similar to the examples plotted in Fig.~\ref{fig:fem_eigen}, the complex eigenvalues in test~(c) are farther away from the real axis than those in the other tests.

In all examples, Algorithm~\ref{code:rand_arnoldi} is able to produce a sufficiently well-con\-di\-tioned Krylov basis $\mat{W}_{m}$, such that the errors of \texttt{rand} and \texttt{rand-ls} are virtually indistinguishable from those of the standard Arnoldi iteration (\texttt{arnoldi}), even with a small sparsity parameter $\zeta$. In contrast, with an incomplete orthogonalization, the Krylov basis $\mat{W}_{m}$ ra\-pi\-dly becomes ill-conditioned, reaching $\kappa_2(\mat{W}_{m}) \geq 10^{15}$ after only $80$ iterations. As a result, \texttt{incomplete} converges significantly more slowly than \texttt{arnoldi}. The basis whitening in \texttt{sFOM} greatly improves the basis conditioning as shown in Table~\ref{tab:conv_diff_full}; however, its condition number is still an order of magnitude higher than the other randomized methods, such that its convergence rate is noticeably slower. We observe that the \texttt{incomplete} method stagnates with $k = 10$ and diverges with $k = 5$ in test~(a). The choice of $k$ only has a minor effect on the convergence of \texttt{sFOM}.

The convergence curves for the restarted methods are shown in Fig.~\ref{fig:conv_diff_restart}. As expected, restarting the Arnoldi procedure slows down the convergence of the method, requiring more iterations to reach the same accuracy. Yet, surprisingly, \texttt{restart-rand} converges faster than the standard \texttt{restart} in test~(a), and even beats the \texttt{deflated} method for $m_r = 20$. This can be explained in terms of the (sketched) Ritz values, as they determine the nodes of the underlying interpolation process for the Krylov approximant $\vec{\hat{f}}_m$ (\ref{eq:rand_arnoldi_approx}).

At restart cycle $k$, let $\Lambda(\mat{R}_{km}) = \{\vartheta_1,\ldots,\vartheta_{km}\}$ and let $q \in \mathcal{P}^{km-1}$ interpolate $f$ at these nodes in the Hermite sense, so that $\vec{\hat{f}}_k = q(\mat{L})\vec{b}$ by~\eqref{eq:rand_arnoldi_polynomial}. Suppose that $D$ is a closed, convex set containing $W(\mat{L})$ and all nodes generated up to cycle $k$, and that $f$ is analytic on a neighborhood of $D$. The Crouzeix--Palencia theorem~\cite{crouzeix_numerical_2017} then gives
\begin{equation}
\label{eq:crouzeix_palencia}
\norm{f(\mat{L})\vec{b} - \vec{\hat{f}}_k} \leq (1 + \sqrt{2}) \; \norm{\vec{b}} \; \max_{z \in W(\mat{L})} |f(z) - q(z)|.
\end{equation}
The interpolation error $|f(z) - q(z)|$ is given by
\begin{equation*}
|f(z) - q(z)| = |f[\vartheta_1, \vartheta_2, ..., \vartheta_{km}; z]| \; \chi_{km}(z)
\end{equation*}
with $\chi_{km}(z) = \Pi_{i = 1}^{km} |z - \vartheta_i|$ and $f[\vartheta_1, \ldots, \vartheta_{km}, z]$ as the divided difference of $f$ of order $km$ over the nodes $\vartheta_1, \ldots, \vartheta_{km}, z$. Using the Hermite-Genocchi representation of the divided difference \cite{palitta_sketched_2025,de_boor_divided_2005} yields
\begin{equation*}
\bigl| f[\vartheta_1,\ldots,\vartheta_{km},z] \bigr| \leq \frac{\max_{\zeta \in D} \bigl| f^{(km)}(\zeta) \bigr|}{(km)!}.
\end{equation*}
Then, substituting the terms in (\ref{eq:crouzeix_palencia}), we have
\begin{equation}
\label{eq:error_bound}
\norm{f(\mat{L})\vec{b} - \vec{\hat{f}}_k} \leq (1 + \sqrt{2}) \; \norm{\vec{b}} \; \frac{\max_{\zeta \in D} \bigl| f^{(km)}(\zeta) \bigr|}{(km)!} \; \max_{z \in W(\mat{L})} |\chi_{km}(z)|.
\end{equation}
From the bound, we can heuristically infer that $\max_{z \in W(\mat{L})} |\chi_{km}(z)|$ will be small when every point of $W(\mat{L})$ is close to one or more eigenvalues of $\mat{R}_{km}$, and that its value is controlled by the most distant point of $W(\mat{L})$ (usually an extreme point) since $\chi_{km}(z)$ grows rapidly outside the convex hull of the nodes. It is worth mentioning that the main purpose of the bound (\ref{eq:error_bound}) is to provide some insight into how the distribution of the interpolation nodes affects the overall convergence of the methods. See \cite{eiermann_analysis_2000,eiermann_restarted_2006,frommer_convergence_2014,palitta_sketched_2025} for a more in-depth analysis.

\begin{figure}[t]
   \centering	\includegraphics[width=0.8\linewidth]{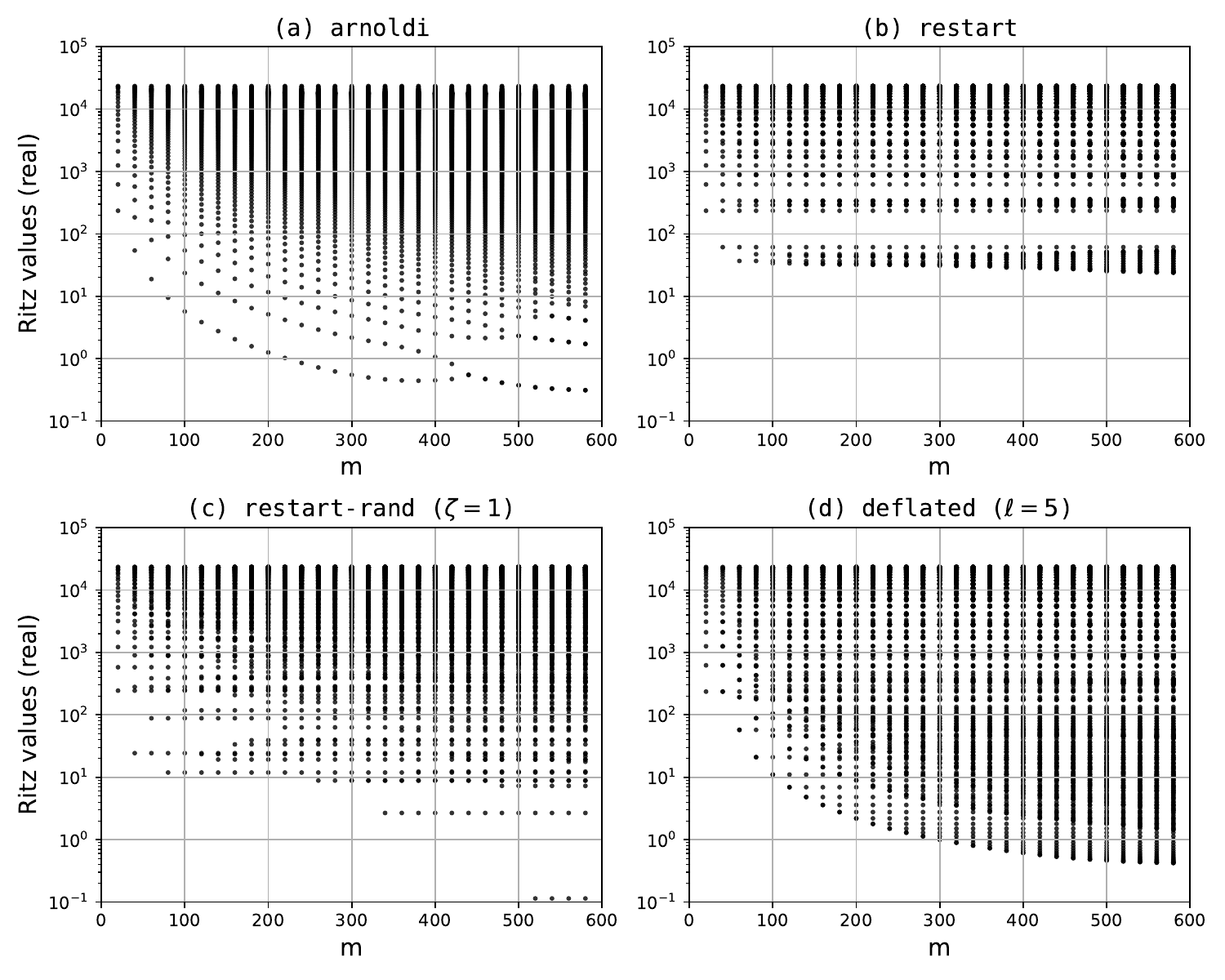}
	\caption{The Ritz values for the convection-diffusion problem with $\alpha = 0.1, \beta = 0.01$, i.e., for test~(a). The restart length was set to $20$.}
	\label{fig:conv_diff_ritz}
\end{figure}

Fig.~\ref{fig:conv_diff_ritz} shows the Ritz values from \texttt{arnoldi}, \texttt{restart}, \texttt{restart-rand}, and \texttt{deflated} for test~(a). Test~(a) has a large value of $\alpha$, such that $\Lambda(\mat{L})$ is quite wide, ranging from $0.120$ to $23,363$ as shown in Table~\ref{tab:conv-diff-eigen}. Most eigenvalues of $\mat{L}$ are real or have a very small imaginary part. After $m = 600$ iterations of the \texttt{arnoldi} method, the Ritz values span the entire spectrum of $\mat{L}$ with the leftmost eigenvalue located at $0.313$ (Fig.~\ref{fig:conv_diff_ritz}a). In contrast, \texttt{restart} has its Ritz values clustered in $20$ discrete regions ranging from $24$ to $23,361$ (see Fig.~\ref{fig:conv_diff_ritz}b), which is the same behavior reported in \cite{afanasjew_generalization_2008,eiermann_restarted_2006,eiermann_deflated_2011}. As a result, the error $|f(z) - q(z)|$ for every point $z \in W(\mat{L})$ that is close to these clusters tends to be smaller than for those that are farther away (e.g., the eigenvalue of $\mat{L}$ located at $0.120$).

From \cite{grigori_restoring_2026,palitta_sketched_2025}, we know that the randomization perturbs the eigenvalues of $\mat{R}_{km}$, such that $W(\mat{R}_{km})$ might be larger than $W(\mat{L})$. Each restart cycle $k$ in \texttt{restart-rand} starts with the perturbed basis vector $\vec{w}^{(k - 1)}_{m + 1}$ from the previous cycle, and thus the perturbations can accumulate over time. Although they are associated with convergence spikes in the unrestarted method \cite{grigori_restoring_2026,palitta_sketched_2025}, they have the opposite effect in the restarted case for this test. With randomization, the eigenvalues of $\mat{R}_{km}$ are now well-distributed over the entire spectrum of $\mat{L}$ with values ranging from $0.114$ to $23,476$. Therefore, the distance from a point in $W(\mat{L})$ to the closest eigenvalue of $\mat{R}_{km}$ is significantly smaller than for \texttt{restart}, especially near the extremes, and thus, the error $|f(z) - q(z)|$ at these points is also lower.

Spectral deflation also breaks the discrete pattern observed for \texttt{restart} by systematically deflating a set of eigenvalues of $\mat{H}_{km}$ at each restart cycle \cite{eiermann_deflated_2011}. Its effectiveness depends on how well the deflated Ritz values approximate the target eigenvalues of $\mat{L}$ at each cycle. As shown in Fig.~\ref{fig:conv_diff_ritz}d, the eigenvalues of $\mat{H}_m$ cover most of the spectrum of $\mat{L}$, with the leftmost one located at $0.423$.

Generally speaking, with a longer restart cycle (i.e., with a large $m_r$), the Ritz values are better distributed over the spectrum of $\mat{L}$, leading to a more accurate approximation of $f(\mat{L}) \vec{b}$. Since $\mat{L}$ has the widest spectrum in test~(a), we see a significant improvement in the convergence rate of the restarted methods after increasing the restart length $m_r$. In the other tests, a good representation of the spectrum of $\mat{L}$ is already attained with a restart length of $20$, and thus increasing $m_r$ has a smaller impact on the convergence of the method compared to test~(a).

\begin{table}[t]
\centering
\caption{Condition number of the basis $\kappa_2(\mat{W}_m)$ and observed distortion $\varepsilon_{obs}$ at the last restart cycle of \texttt{restart-rand} for different combinations of $d$ and $\zeta$. The restart length $m_r$ was set to $20$.}
\label{tab:conv_diff_restart}
\begin{tabular}{cccccccc}
\hline
 &  & \multicolumn{2}{c}{(a) $\alpha = 0.1, \beta = 0.01$} & \multicolumn{2}{c}{(b) $\alpha = 0.01, \beta = 0.01$} & \multicolumn{2}{c}{(c) $\alpha = 0.01, \beta = 1$} \\
\cline{3-4}\cline{5-6}\cline{7-8}
$d$ & $\zeta$ & $\kappa_2(\mat{W}_m)$ & $\varepsilon_{obs}$ & $\kappa_2(\mat{W}_m)$ & $\varepsilon_{obs}$ & $\kappa_2(\mat{W}_m)$ & $\varepsilon_{obs}$ \\
\hline
$160$ & $1$ & $6.28 \pm 5.95$ & $0.90 \pm 0.18$ & $6.77 \pm 5.76$ & $0.95 \pm 0.24$ & $6.49 \pm 5.26$ & $0.94 \pm 0.23$ \\
$320$ & $1$ & $1.52 \pm 0.04$ & $0.44 \pm 0.03$ & $1.60 \pm 0.04$ & $0.48 \pm 0.06$ & $2.78 \pm 2.69$ & $0.59 \pm 0.24$ \\
$640$ & $1$ & $1.36 \pm 0.03$ & $0.33 \pm 0.02$ & $1.40 \pm 0.08$ & $0.40 \pm 0.12$ & $1.34 \pm 0.01$ & $0.30 \pm 0.02$ \\
$1280$ & $1$ & $1.31 \pm 0.13$ & $0.27 \pm 0.10$ & $1.27 \pm 0.03$ & $0.25 \pm 0.03$ & $1.25 \pm 0.01$ & $0.22 \pm 0.01$ \\
$2560$ & $1$ & $1.17 \pm 0.01$ & $0.17 \pm 0.01$ & $1.20 \pm 0.04$ & $0.19 \pm 0.03$ & $1.17 \pm 0.01$ & $0.16 \pm 0.01$ \\
$5120$ & $1$ & $1.12 \pm 0.02$ & $0.12 \pm 0.02$ & $1.13 \pm 0.03$ & $0.13 \pm 0.03$ & $1.11 \pm 0.01$ & $0.12 \pm 0.01$ \\
\hline
$320$ & $1$ & $1.52 \pm 0.04$ & $0.44 \pm 0.03$ & $1.60 \pm 0.04$ & $0.48 \pm 0.06$ & $2.78 \pm 2.69$ & $0.59 \pm 0.24$ \\
$320$ & $4$ & $1.60 \pm 0.03$ & $0.50 \pm 0.03$ & $1.56 \pm 0.04$ & $0.47 \pm 0.04$ & $1.61 \pm 0.07$ & $0.54 \pm 0.05$ \\
$320$ & $8$ & $1.57 \pm 0.03$ & $0.50 \pm 0.06$ & $1.57 \pm 0.04$ & $0.49 \pm 0.05$ & $1.58 \pm 0.08$ & $0.50 \pm 0.07$ \\
$320$ & $16$ & $1.59 \pm 0.03$ & $0.50 \pm 0.04$ & $1.55 \pm 0.07$ & $0.47 \pm 0.07$ & $1.60 \pm 0.08$ & $0.50 \pm 0.04$ \\
$320$ & $32$ & $1.54 \pm 0.07$ & $0.47 \pm 0.05$ & $1.60 \pm 0.05$ & $0.50 \pm 0.06$ & $1.58 \pm 0.03$ & $0.48 \pm 0.04$ \\
\hline
\end{tabular}
\end{table}

\begin{figure}[t]
	\centering	\includegraphics[width=\textwidth]{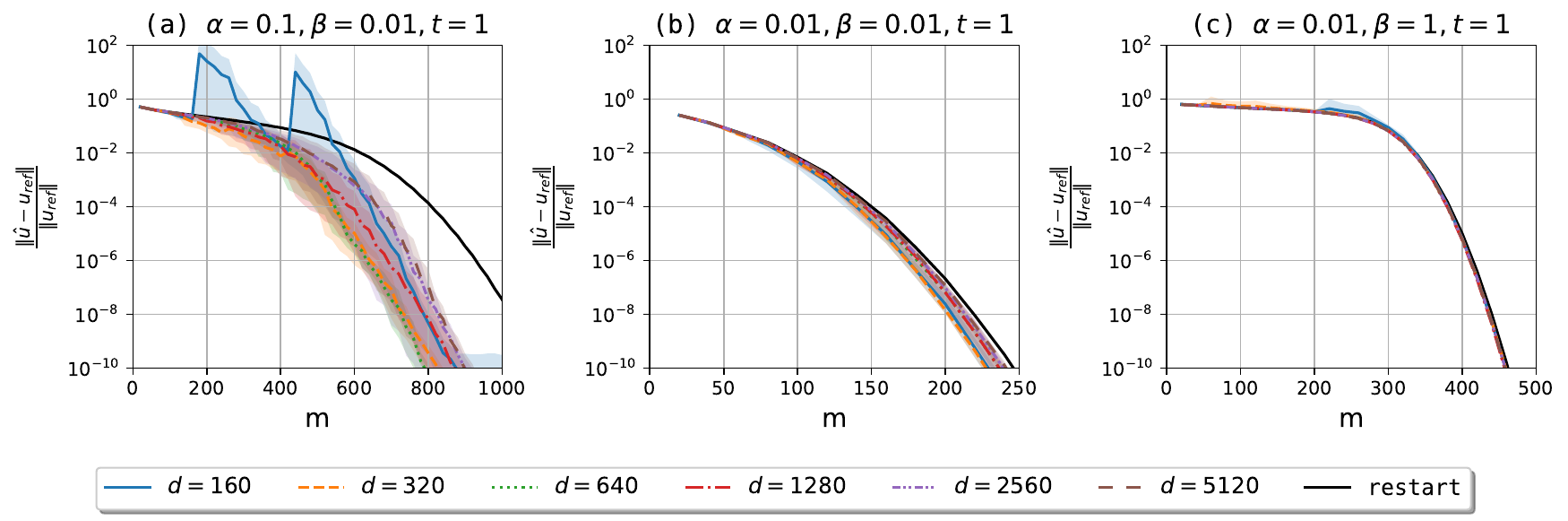}
	\caption{Convergence curves for \texttt{restart-rand} for different sketching dimensions $d$ for the convection-diffusion example. The restart length $m_r$ was set to $20$ and the sparsity parameter $\zeta$ was set to $1$.}
	\label{fig:conv_diff_restart_s}
\end{figure}

Fig.~\ref{fig:conv_diff_restart_s} shows the effect of the sketching dimension $d$ on the convergence of \texttt{restart-rand}. The basis conditioning $\kappa_2(\mat{W}_m)$ and the observed distortion $\varepsilon_{obs}$ at the last restart cycle are reported in Table~\ref{tab:conv_diff_restart}. With $d = 160$, the random sketching introduces such a high distortion $\varepsilon$ that it destabilizes the restarted method and causes the convergence spikes as shown in Fig.~\ref{fig:conv_diff_restart_s}a, although it still manages to converge to the reference solution with an error of $\sim 10^{-13}$. Despite the high distortion, the basis $\mat{W}_m$ is relatively well conditioned with $\kappa_2(\mat{W}_m) \leq 12$, although this is considerably higher than for the other sketching dimensions. The convergence spikes disappear for $d \geq 320$ since the distortion is now around or below the recommended value $\varepsilon \leq 0.5$ \cite{balabanov_randomized_2022}. Higher values of $d$ seem to move the curve towards the regular restarted method in tests~(a) and (b). The sparsity parameter $\zeta$ does not affect the convergence of the method in any of our experiments. The basis conditioning and distortion are also fairly similar across tests~(a) and (b). In test~(c), $\zeta = 1$ produces a larger mean $\kappa_2(\mat{W}_m)$ and $\varepsilon_{obs}$ than the others, but the convergence rate does not seem to be negatively affected in our tests.

\subsection{Circular Membrane}
\label{subsec:wave}

The second example consists of simulating the vibrations of a circular membrane \cite[Chapter 9]{strutt_theory_2011}. Let us consider a membrane $\Omega$ of radius $\mu$ centered at the origin. At a time $t > 0$, the height of the membrane at any point $(x, y)$ is given by $u(x, y, t)$, measured from the rest position. The membrane is attached to a rigid frame, such that $u(x, y, t) = 0$ for all $(x, y) \in \partial \Omega$. The vibrations in the membrane can then be described as
\begin{equation}
\label{eq:drum_def}
\frac{\partial^2}{\partial t^2} u(x, y, t) = \nu^2 \, \Delta u(x, y, t),
\end{equation}
where $\nu > 0$ is the speed at which the transverse waves propagate in the membrane.
After using the finite element method \cite{zienkiewicz_finite_2013} to describe (\ref{eq:drum_def}) in terms of a discrete mesh over $\Omega$, we have
\begin{equation}
\label{eq:drum_discrete}
\frac{d^2}{d t^2} \vec{\hat{u}} = -\nu^2 \,\mat{L} \vec{\hat{u}},
\end{equation}
where $\mat{L} = \mat{M}^{-1} \mat{A}$ for the stiffness matrix $\mat{A}$ and the lumped mass matrix $\mat{M}$. Suppose that the initial conditions are set to
\begin{equation}
\label{eq:drum_initial}
u(x, y, 0) = J_p \left (\frac{\eta_{pk} \sqrt{x^2 + y^2}}{\mu} \right), \qquad \frac{d}{dt}u(x, y, 0) = 0,
\end{equation}
where $J_p (x)$ is the Bessel function
\begin{equation*}
J_p(x) = \sum_{k = 0}^\infty{\frac{(-1)^k}{k! \, (k + p)!} \left(\frac{x}{2} \right)^{2k + p}},
\end{equation*}
and $\eta_{pk}$ is the $k$-th zero of $J_p (x)$. Then, the solution for (\ref{eq:drum_discrete}) can be written as
\begin{equation}
\label{eq:drum_sol}
\vec{\hat{u}}(t) = \cos(\nu t \, \sqrt{\mat{L}}) \, \vec{b},
\end{equation}
where $\sqrt{\mat{L}}$ is any square root of $\mat{L}$ \cite[p. 124]{gantmakher_theory_1959} and $\vec{b}$ is a vector containing the value of $J_p (x)$ for each node in the discrete mesh. For our experiments, we consider a circular membrane with radius $\mu = 1$ and set the initial conditions using a fourth-order Bessel function $J_4 (x)$ and its fourth zero. The matrix $\mat{L}$ was generated with \texttt{gmsh} \cite{geuzaine_gmsh_2009} and \texttt{FreeFem++}~\cite{hecht_new_2012} using P1 finite elements with a mesh size of $0.0025$. To impose the Dirichlet boundary conditions $u(x, y, t) = 0$ for all $(x, y) \in \partial \Omega$, we replace the corresponding rows of $\mat{L}$ with the identity matrix \cite{larson_finite_2013}. The resulting matrix $\mat{L}$ has $582,547$ rows and $4,062,636$ nonzeros. The maximum and minimum eigenvalues of $\mat{L}$ are $4.502 \times 10^6$ and $1$, respectively. Its condition number $\kappa_2(\mat{L})$ is $7.11 \times 10^7$. We adopted as reference the solution obtained by \texttt{arnoldi} with $m = 2000$. The matrix cosine was computed using the scaling-and-squaring method from \cite{al-mohy_new_2015}.

\begin{figure}[t]
\begin{minipage}[t]{0.49 \linewidth}
    \centering	\includegraphics[width=\textwidth]{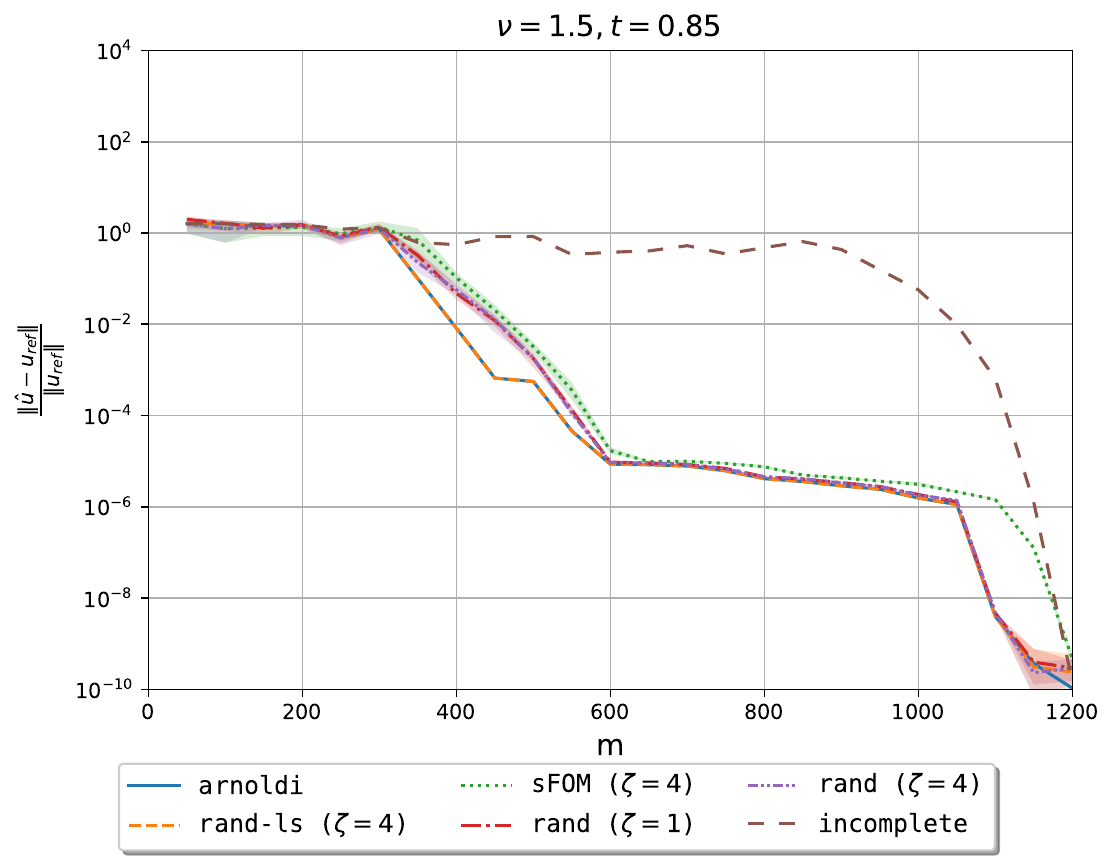}
	\caption{Convergence curves for non-restarted methods for the circular membrane example.
	The truncation parameter $k$ was set to $80$ and the sketching dimension $d$ was set to $3600$.}
	\label{fig:drum_convergence_full}
\end{minipage}
\hfill
\begin{minipage}[t]{0.49 \linewidth}
    \centering	\includegraphics[width=\textwidth]{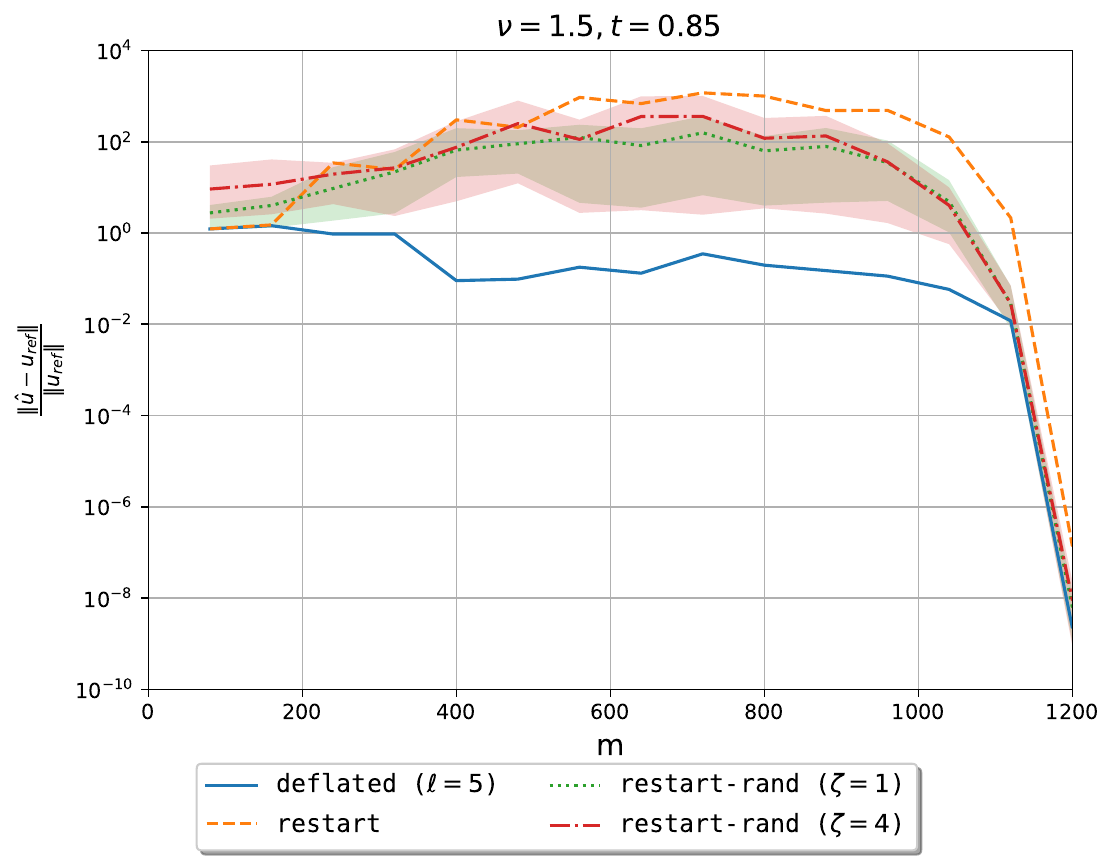}
	\caption{Convergence curves for restarted methods when simulating the circular membrane example. The restart length $m_r$ was set to $80$ and the sketching dimension $d$ was set to $960$.}
	\label{fig:drum_convergence_restart}
\end{minipage} \\
\end{figure}

\begin{table}[t]
\centering
\caption{The relative $\ell_2$ error, condition number of the basis $\kappa_2(\mat{W}_m)$, and observed distortion $\varepsilon_{obs}$ at $m = 1200$. The restart length $m_r$ was set to $80$. For \texttt{sFOM}, we report $\kappa_2(\mat{W}_m \mat{T}_m^{-1})$ instead of $\kappa_2(\mat{W}_m)$.}
\label{tab:drum-sketching}
\begin{tabular}{cccccc}
\hline
method & $d$ & $\zeta$ & $\frac{\norm{\hat{\vec{u}} - \vec{u}_{ref}}}{\norm{\vec{u}_{ref}}}$ & $\kappa_2(\mat{W}_m)$ & $\varepsilon_{obs}$ \\
\hline
\texttt{rand} & $3600$ & $1$ & $(2.97 \pm 1.34) \times 10^{-10}$ & $4.61 \pm 0.17$ & $0.91 \pm 0.01$ \\
\texttt{rand} & $3600$ & $4$ & $(2.93 \pm 1.35) \times 10^{-10}$ & $3.73 \pm 0.04$ & $0.87 \pm 0.00$ \\
\texttt{rand-ls} & $3600$ & $4$ & $(2.40 \pm 2.08) \times 10^{-10}$ & $3.73 \pm 0.04$ & $0.87 \pm 0.00$ \\
\texttt{sFOM} & $3600$ & $4$ & $(4.79 \pm 1.15) \times 10^{-10}$ & $17.52 \pm 3.44$ & $0.99 \pm 0.00$ \\
\hline
\texttt{restart-rand} & $160$ & $1$ & $(3.71 \pm 5.12) \times 10^{9}$ & $5.92 \pm 0.40$ & $0.94 \pm 0.01$ \\
\texttt{restart-rand} & $240$ & $1$ & $(8.46 \pm 18.91) \times 10^{4}$ & $3.85 \pm 0.19$ & $0.87 \pm 0.01$ \\
\texttt{restart-rand} & $480$ & $1$ & $(3.44 \pm 7.50) \times 10^{-6}$ & $2.38 \pm 0.17$ & $0.70 \pm 0.03$ \\
\texttt{restart-rand} & $960$ & $1$ & $(8.75 \pm 4.68) \times 10^{-9}$ & $1.91 \pm 0.21$ & $0.56 \pm 0.07$ \\
\texttt{restart-rand} & $1920$ & $1$ & $(2.31 \pm 3.32) \times 10^{-8}$ & $1.56 \pm 0.12$ & $0.41 \pm 0.06$ \\
\texttt{restart-rand} & $3840$ & $1$ & $(1.38 \pm 1.56) \times 10^{-8}$ & $1.36 \pm 0.04$ & $0.30 \pm 0.02$ \\
\texttt{restart-rand} & $7680$ & $1$ & $(2.67 \pm 3.37) \times 10^{-8}$ & $1.23 \pm 0.04$ & $0.20 \pm 0.03$ \\
\hline
\texttt{restart-rand} & $960$ & $4$ & $(9.43 \pm 11.20) \times 10^{-9}$ & $1.77 \pm 0.03$ & $0.51 \pm 0.01$ \\
\texttt{restart-rand} & $960$ & $8$ & $(2.66 \pm 2.24) \times 10^{-8}$ & $1.77 \pm 0.02$ & $0.51 \pm 0.01$ \\
\texttt{restart-rand} & $960$ & $16$ & $(7.35 \pm 14.11) \times 10^{-8}$ & $1.77 \pm 0.02$ & $0.52 \pm 0.01$ \\
\texttt{restart-rand} & $960$ & $32$ & $(1.35 \pm 0.89) \times 10^{-8}$ & $1.77 \pm 0.03$ & $0.52 \pm 0.01$ \\
\hline
\end{tabular}
\end{table}

Figs.~\ref{fig:drum_convergence_full} and \ref{fig:drum_convergence_restart} show the convergence curves for all numerical methods. Table~\ref{tab:drum-sketching} reports the error, basis conditioning, and observed distortion at $m = 1200$. As a consequence of the oscillatory nature of the problem, the classical \texttt{arnoldi} method has a staircase-shaped convergence and eventually stagnates around $10^{-10}$. The errors of \texttt{rand} and \texttt{rand-ls} closely follow those of the classical method, while \texttt{incomplete} takes around $m = 900$ iterations to start converging to the solution due to a very ill-conditioned basis ($\kappa_2(\mat{W}_m) \geq 10^{16}$). It achieves the same accuracy as the other methods at the end of the experiment. The \texttt{sFOM} curve is similar to those of the other randomized methods, with a slight convergence delay around the $m = 1000$ mark.

The errors of \texttt{restart} and \texttt{restart-rand} are fairly large in the first $12$ restart cycles (i.e., $m \leq 960$), but they quickly converge towards the solution thereafter. The numerical instability is typically caused by high-frequency oscillations. Spectral deflation is able to prevent the method from diverging; however, it still stagnates around $10^{-1}$ for $12$ restart cycles before converging to the solution. In this experiment, the restart length controls how far the method diverges from the reference solution; however, it does not change the overall convergence rate of the method.
Regarding the sketching dimension $d$, the observed distortion $\varepsilon_{obs}$ is quite large for $d = 160$ and $d = 240$, which causes the method to diverge, even though the basis is relatively well conditioned ($\kappa_2(\mat{W}_m) \leq 6$). The method converges for $d = 480$, yet it can only produce a solution within $10^{-6}$ of the reference regardless of the number of restarts. Any sketching dimension $d \geq 960$ has a sufficiently low distortion that \texttt{restart-rand} has an error of $10^{-8}$ on the last cycle. The sparsity parameter $\zeta$ does not seem to affect the overall convergence of the method in our tests.

\subsection{Graph Laplacian}
\label{subsec:laplacian}

\afterpage {

\begin{table}[H]
\centering
\caption{Properties of complex networks in the graph Laplacian example. Here, M stands for millions.}
\label{tab:graphs}
\begin{tabular}{lllllp{0.4\textwidth}}
\toprule
Graph              & Nodes        & Edges         & Is directed? & $\lambda_{max}(\mat{L})$ & Description                                                                                                      \\ \midrule
\texttt{kronecker} & $8.9$M  & $529$M & No           & $406,068$                & Kronecker graph used by the Graph500 benchmark \cite{graph500,leskovec_kronecker_2010}                           \\ \midrule
\texttt{orkut}     & $3.1$M  & $237$M & No           & $33,314$                 & Social network of Orkut users in 2007 \cite{leskovec_snap,mislove_measurement_2007}                              \\ \midrule
\texttt{uk-2002}   & $18.5$M & $316$M & Yes          & $2,449$                  & Web graph of the \textit{.uk} domain in 2002 \cite{boldi_ubicrawler_2004,boldi_layered_2011,boldi_webgraph_2004} \\
\bottomrule
\end{tabular}
\end{table}
\vfill
\begin{figure}[H]
\begin{minipage}[t]{\linewidth}
    \centering	\includegraphics[width=\textwidth]{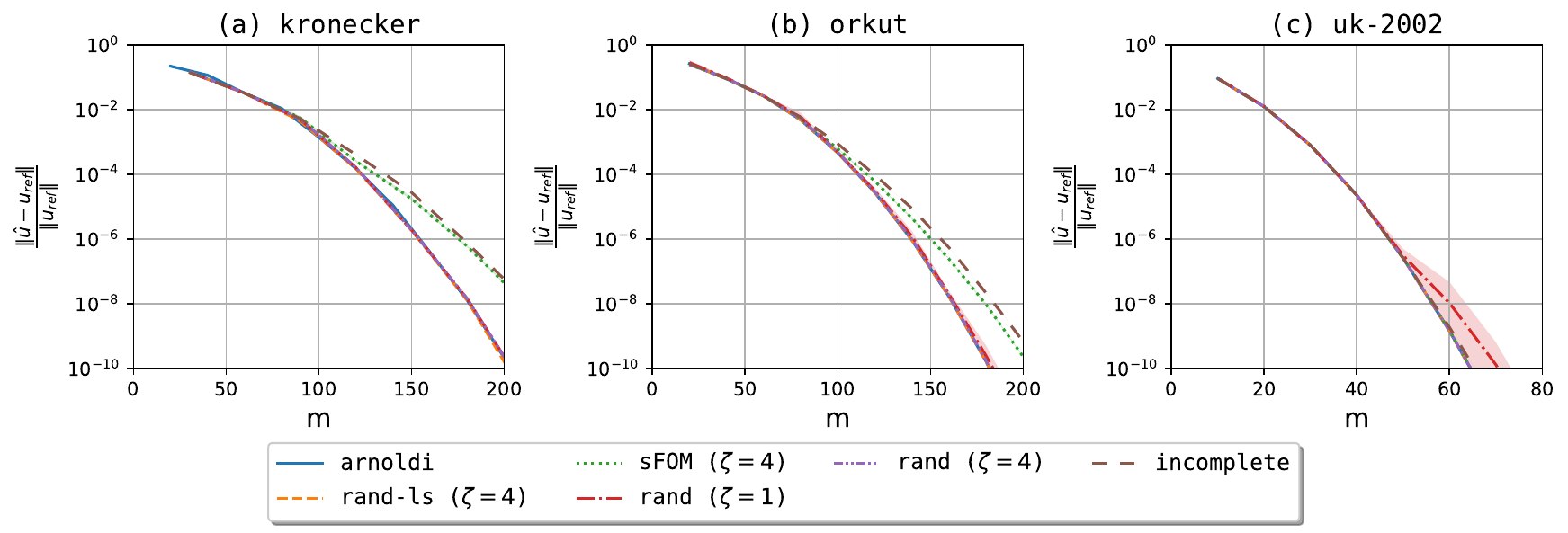}
	\caption{Convergence curves for non-restarted methods for the graph Laplacian example with $t = 0.1$.
	The truncation parameter $k$ was set to $40$ and the sketching dimension $d$ was set to $1200$.}
	\label{fig:laplacian_convergence_full}
\end{minipage}
\vfill
\begin{minipage}[t]{\linewidth}
    \centering	\includegraphics[width=\textwidth]{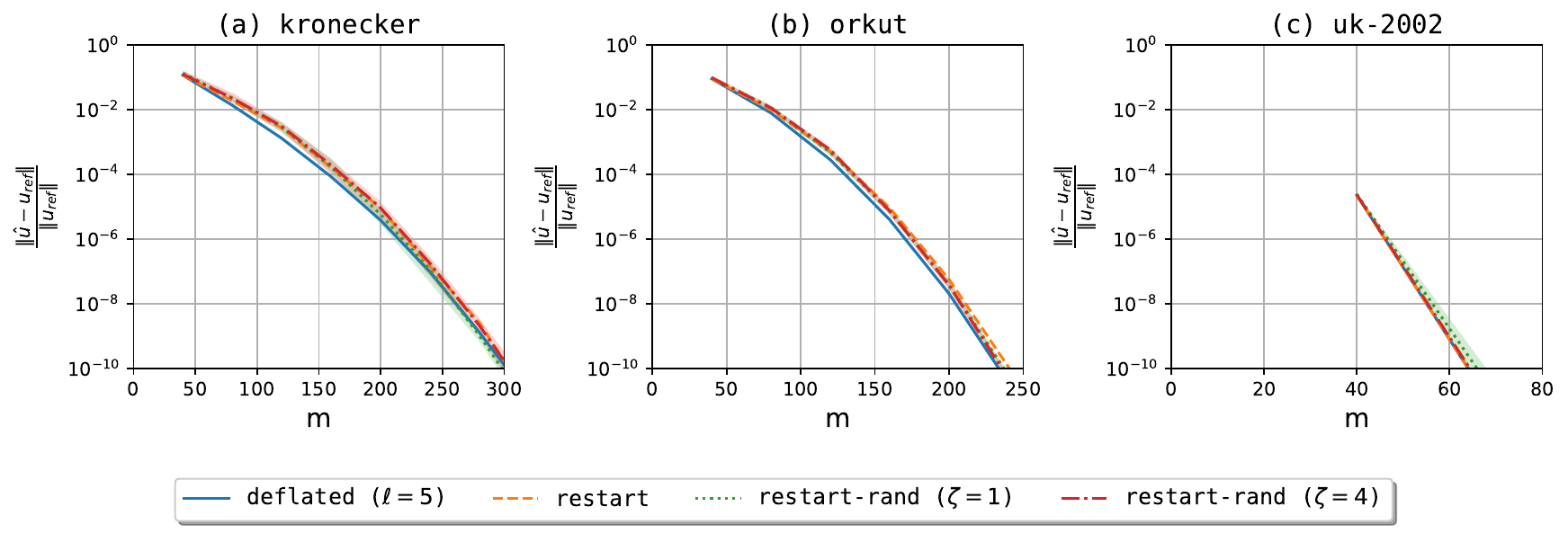}
	\caption{Convergence curves for restarted methods for the graph Laplacian example with $t = 0.1$. The restart length $m_r$ was set to $40$ and the sketching dimension $d$ was set to $320$.}
	\label{fig:laplacian_convergence_restart}
\end{minipage} \\
\end{figure}
}

The third and last example consists of modeling the diffusion on a graph $G = (V, E)$ as
\begin{equation}
\label{eq:heat_network}
\frac{d \vec{x}}{dt} = -\mat{L} \vec{x}, \quad \vec{x}(0) = \vec{x}_0,
\end{equation}
where $\mat{L}$ is the out-degree graph Laplacian (i.e., $\mat{L}[i, i] = \sum_{i \neq j} |\mat{L}[i, j]|$ and $\mat{L}[i, j] = -1$ if there is an edge from node $i$ to node $j$, $0$ otherwise) and $\vec{x}$ is the probability distribution over the nodes at a time $t \geq 0$. Note that the Laplacian is singular ($\mat{L} \vec{1} = \vec{0}$) and its entire spectrum is located in the closed right half-plane. The solution for (\ref{eq:heat_network}) is then expressed as
\begin{equation}
\label{background:eq:heat_laplacian_sol}
\vec{x} = e^{-t \mat{L}} \vec{x}_0 = \mat{H}(t) \vec{x}_0.
\end{equation}
The matrix-valued function $\mat{H}(t)$ is called the \textit{heat kernel} on $G$ and models the diffusion over the network (see \cite[Section 6.2]{benzi_matrix_2020}). In this example, the initial probability distribution $\vec{x}_0$ was set by first drawing $n$ i.i.d.\ random numbers, collecting them into a vector $\vec{y} = \{ y_i \sim \textsc{uniform}(0, 1), \quad i = 1, 2, ..., n \}$, and then taking $\vec{x}_0 = \vec{y} / \norm{\vec{y}}_1$. A summary of the complex networks used in this section is given in Table~\ref{tab:graphs}. The reference solution was obtained by running the \texttt{restart} method with a restart length $m_r = 100$ and a tolerance of $10^{-12}$. Note that most results from this section can be extended to the in-degree case with some minor modifications.

Figs.~\ref{fig:laplacian_convergence_full} and \ref{fig:laplacian_convergence_restart} show the convergence curves for all numerical methods. For both the \texttt{kronecker} and \texttt{orkut} networks, the error curves are dictated by the type of orthogonalization --- (sketched) full orthogonalization (\ie, \texttt{arnoldi} and \texttt{rand}/\texttt{rand-ls}), incomplete orthogonalization, and restarting --- without a clear separation between the standard and randomized versions. For \texttt{uk-2002}, all the methods share the same convergence rate due to the narrow spectrum of $\mat{L}$.

\subsection{Performance}
\label{sec:performance}

\begin{figure}[t]
	\centering	\includegraphics[width=\textwidth]{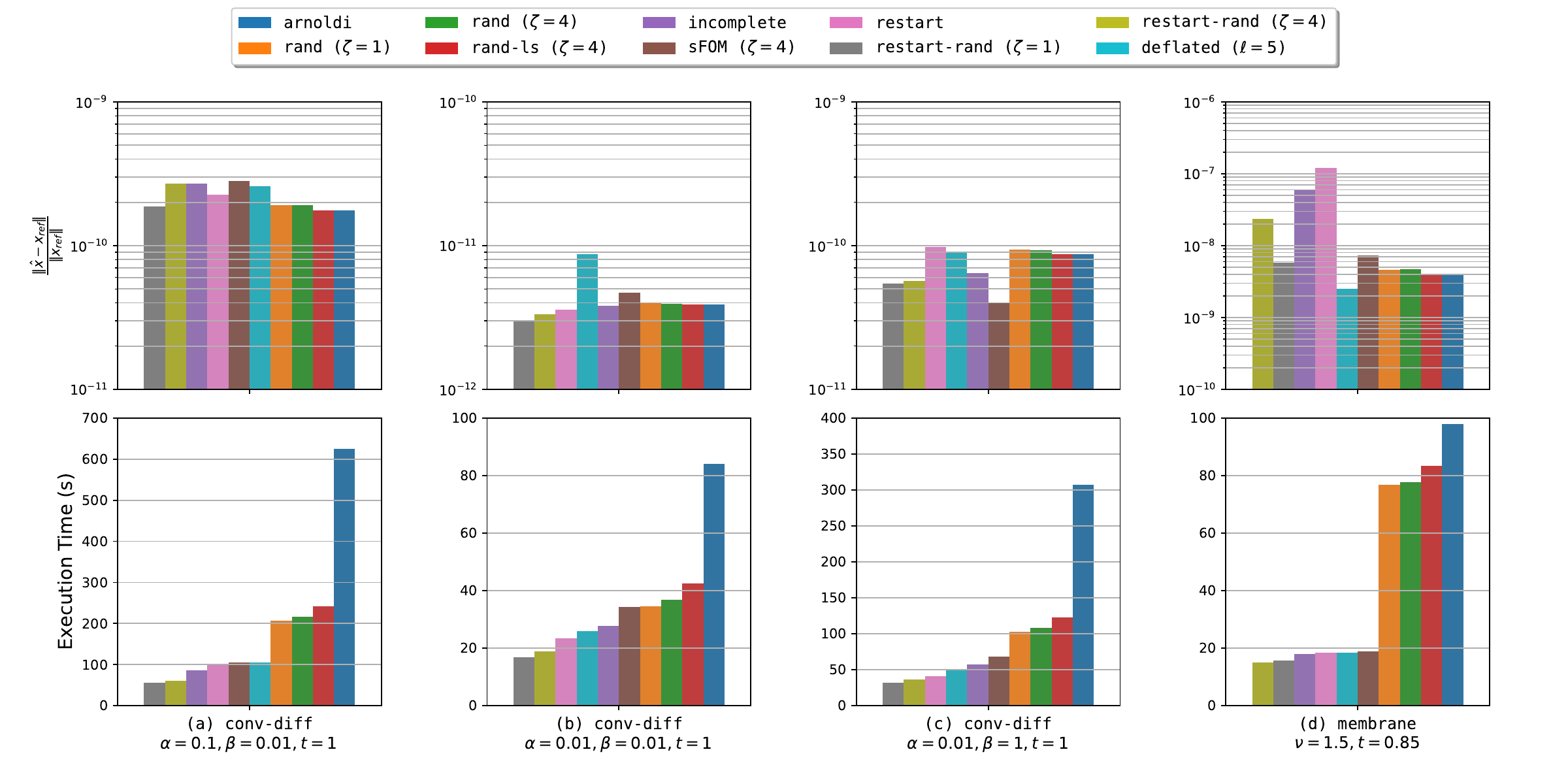}
	\caption{Comparison between the execution time and accuracy of the different Krylov methods.}
	\label{fig:timing}
\end{figure}

The performance tests were carried out on a commodity desktop with an AMD Ryzen 7 5800X 8C/16T @4.6GHz and 32GB of RAM, running Arch Linux. Fig.~\ref{fig:timing} compares the execution time and accuracy of the numerical methods. For the \texttt{conv-diff} test~(a--c), we set $m_r = k = 20$; $d = 3200$ for \texttt{rand}, \texttt{rand-ls}, and \texttt{sFOM}; and $d = 320$ for \texttt{restart-rand}. For the \texttt{membrane} test~(d), we set $m_r = k = 80$; $d = 3600$ for \texttt{rand}, \texttt{rand-ls}, and \texttt{sFOM}; and $d = 960$ for \texttt{restart-rand}. The number of restarts and the basis size $m$ were adjusted in such a way that the accuracy of all methods is as close as possible.

In all tests, \texttt{arnoldi} has the highest execution time among all methods due to the full orthogonalization of the Krylov basis. Replacing the classical Arnoldi procedure with the randomized version (i.e., the \texttt{rand} method) leads to a speedup of $1.26 - 2.90\times$ for $\zeta = 4$ and $1.28 - 3.03\times$ for $\zeta = 1$. Generally speaking, these methods are quite memory-bound, and thus, the performance gains from \texttt{rand} heavily depend on the memory subsystem. The \texttt{conv-diff} examples are quite large, i.e., $\mat{L}$ consumes $1.13$GB, while a single column of $\mat{V}_m$/$\mat{W}_m$ uses $38.4$MB. This entails that the method needs to load these matrices from the main memory as the L3 cache has only $32$MB. The modified Gram-Schmidt process orthogonalizes the basis via a sequence of BLAS-1 routines (\texttt{axpy} and \texttt{dot}), while RGS updates the basis $\mat{W}_m$ using a single BLAS-2 routine (\texttt{gemv}), which is significantly more memory efficient. The sketch $\mat{S W}_m$ is sufficiently small to fit into the L3 cache, e.g., a $3200 \times 600$ matrix $\mat{S W}_m$ consumes just $15.4$MB, and a few of its columns can be stored comfortably in L2. As a result, \texttt{rand} has a speedup of $2\times$ or more for this test. In contrast, \texttt{membrane} is noticeably smaller with $\mat{L}$ and $\vec{v}_k$/$\vec{w}_k$ using $65$MB and $4.66$MB, respectively. Hence, calling \texttt{axpy} after a \texttt{dot} product is relatively cheap since the working vectors are already loaded in the cache, and thus, the performance gains of \texttt{rand} over \texttt{arnoldi} are smaller ($26-28\%$) for this test. Solving the least-squares problem using LSMR for \texttt{rand-ls} imposes a $9$ to $23\%$ performance penalty depending on the size of the input matrix and Krylov basis.

With a restarted Krylov method, the program spends significantly less time on the orthogonalization of the basis, since it only needs to work with a small set of basis vectors at each restart cycle. As a result, the faster orthogonalization in Algorithm~\ref{code:rand_arnoldi} has less impact on the overall performance of the program. For tests~(b), (c), and (d), where all restarted methods have similar convergence rates, \texttt{restart-rand} with $\zeta = 4$ reduces the execution time by $13-25\%$ relative to the standard \texttt{restart}. Reducing the sparsity parameter to $\zeta = 1$, it now achieves a speedup of $1.18\times$ to $1.40\times$. Spectral deflation actually hurts the performance in these tests by up to $20\%$ since it requires computing and sorting the Ritz values, augmenting the subspace with Schur vectors, and orthogonalizing $\ell = 5$ additional vectors. It is worth mentioning that the relative error of the regular \texttt{restart} method did not fall below $10^{-7}$ in the \texttt{membrane} test. Other methods do not have this problem. In test~(a), \texttt{restart-rand} shows a speedup of $1.61$ and $1.77$ over the standard \texttt{restart} method with $\zeta = 4$ and $\zeta = 1$, respectively. Recall from Section~\ref{subsec:conv_diff} that the convergence rate of \texttt{restart-rand} is significantly higher than \texttt{restart} in test~(a), and thus requires fewer restarts to achieve the same accuracy. Despite having a fast convergence rate, \texttt{deflated} has similar performance to the regular \texttt{restart}.

Both \texttt{incomplete} and \texttt{sFOM} have a fixed orthogonalization cost per new basis vector, regardless of the number of vectors already generated, but require a larger number of iterations to achieve the target accuracy compared to \texttt{arnoldi}. As a result, their performance is similar to the restarted procedures. \texttt{sFOM} pays a small overhead for sketching and basis whitening, increasing the execution time by between $4\%$ and $19\%$.

\begin{figure}[t]
	\centering	\includegraphics[width=\textwidth]{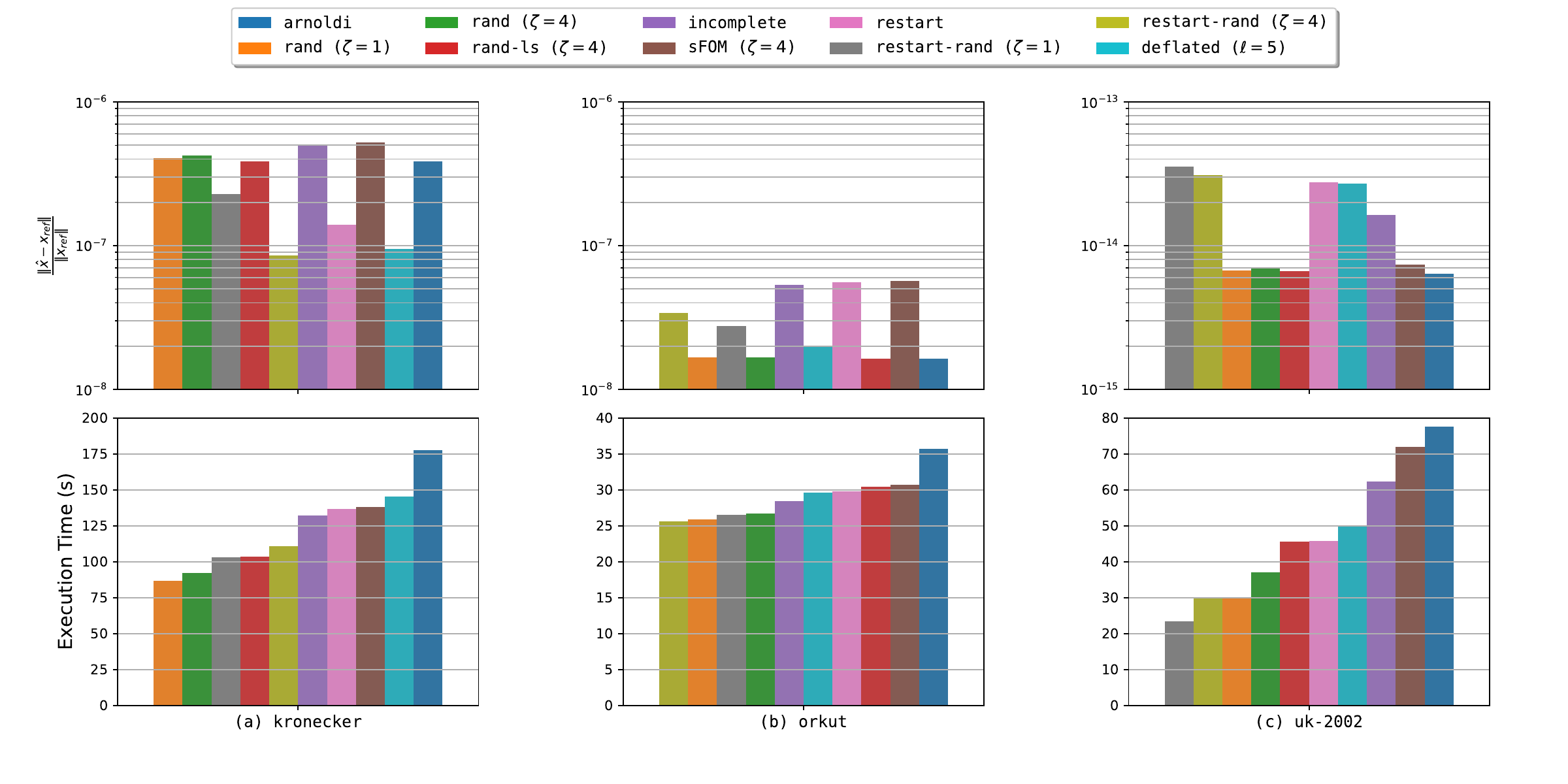}
	\caption{Comparison between the execution time and accuracy of the different Krylov methods for the graph Laplacian example. For this test, we set $m_r = k = 40$; $d = 1200$ for \texttt{rand}, \texttt{rand-ls}, and \texttt{sFOM}; and $d = 320$ for \texttt{restart-rand}.}
	\label{fig:timing_laplacian}
\end{figure}

Fig.~\ref{fig:timing_laplacian} shows the execution time of all methods for solving the graph Laplacian example. For the \texttt{kronecker} and \texttt{orkut} networks, the Laplacian matrix $\mat{L}$ is noticeably larger and/or denser, such that constructing the next basis vector $\vec{v}_{k + 1} = \mat{L} \vec{v}_k$ is significantly more expensive. As restarting delays the convergence of the method, it requires additional sparse matrix-vector products with $\mat{L}$; hence, the benefit of a faster orthogonalization disappears. It is a similar story for \texttt{sFOM} and \texttt{incomplete}. In this example, the fastest method is actually the unrestarted randomized version, \texttt{rand} with $\zeta = 1$. Nevertheless, a restarted method is still beneficial in terms of memory consumption: \texttt{restart} consumes around $11.1$GB of memory for the \texttt{kronecker} network, while \texttt{arnoldi} requires $18.8$GB. The gap becomes wider for a larger number of iterations and/or a larger matrix. In comparison, all numerical methods share the same convergence rate for \texttt{uk-2002}, requiring only $80$ iterations (or $2$ restart cycles) to reach an error of less than $10^{-12}$ (the tolerance of the reference solution). Therefore, the main differentiator here is the orthogonalization cost. \texttt{restart-rand} with $\zeta = 1$ has the lowest execution time ($23$s) followed by \texttt{restart-rand} with $\zeta = 4$ and \texttt{rand} with $\zeta = 1$ ($30$s).

\section{Conclusion}
\label{sec:conclusion}
In this article, we combine the restart technique with the randomized Krylov method from \cite{timsit_randomized_2023} in the context of general matrix functions. Our numerical experiments show that our randomized algorithm significantly outperforms the classical method while producing a highly accurate solution for complex problems. In some cases, the randomization can even lead to better convergence rates. Although our method performs well in practice, many of its features are not completely understood. In particular, it is missing a formal explanation of how the Ritz value perturbations can improve the convergence rate of restarted methods under certain conditions. Likewise, a theoretical analysis of its convergence may lead to a better error estimation, allowing greater control over the number of restarts. Our randomized Krylov method can also be combined with other acceleration techniques, such as spectral deflation \cite{eiermann_deflated_2011}, weighted inner products \cite{embree_weighted_2017}, or ``thick'' restarting \cite{baker_technique_2005}, to further improve its performance.

\begin{appendix}
\section{Implementation details}
\label{appendix:implementation}
\setcounter{equation}{0}
\renewcommand{\theequation}{A.\arabic{equation}}

The test suite was written using the HighFM toolbox\footnote{\url{https://gitlab.com/highfm/highfm}}. This toolbox is written entirely in \texttt{C++20}, and its main goal is to provide a user-friendly interface for solving problems from linear algebra without sacrificing performance. More specifically:

\textbf{Matrices and their operations:} Sparse matrices follow the compressed sparse row (CSR) format, while dense matrices are stored in column-major form. Matrix and vector operations are directly translated into routine calls to Intel MKL with some operations hand-optimized to achieve the highest performance on the target system.

\textbf{Random sketching:} We use \texttt{PCG64DXSM} \cite{oneill_pcg_2014} as the pseudorandom number generator. It is an improved version of the default generator in the popular NumPy \cite{numpy} module. The sparse sign matrix is generated at the start of each randomized method and stored in a compressed sparse row format. It is then applied to the target vector or matrix using \texttt{spmv} and \texttt{spmm}, respectively.

\textbf{Matrix functions:} $\varphi_1(\mat{A})$, $\exp(\mat{A})$, and $\cos(\sqrt{\mat{A}})$ were all computed using the scaling-and-squaring method with the Padé approximant from \cite{higham_scaling_2005,al-mohy_new_2015}. The degree of the approximant is automatically chosen through a careful analysis of the error bounds. The LAPACK routines used are \texttt{dsgesv} for LU solve and \texttt{dgees} for the Schur decomposition.

\textbf{sFOM:} For \texttt{sFOM}, the program calls \texttt{dgeqrf} for the thin QR decomposition \linebreak ($\mat{SW}_m = \mat{Q}_m \mat{T}_m$) and \texttt{dtrsm} for the triangular linear solve. Per LAPACK specification, the factors of the QR decomposition are stored in a single object. The original formula \cite{guttel_randomized_2023} for \texttt{sFOM} is $\vec{f}_m^{gs} = \mat{W}_m \mat{T}_m^{-1} \, f(\mat{Q}_m^* \mat{S AW}_m \mat{T}_m^{-1}) \, \mat{Q}_m^*  \mat{S} \vec{b}$. We use the alternative form proposed in \cite[Eq. 18]{palitta_sketched_2025}:
\begin{equation}
\vec{f}_m^{gs} = \| \mat{S} \vec{b} \| \mat{W}_m \mat{T}_m^{-1} \, f(\mat{X}_m) \, \vec{e}_1, \qquad \mat{X}_m = \mat{T}_m \mat{R}_m \mat{T}_m^{-1} + \frac{r_{m + 1, m}}{t_m} \vec{t} \vec{e}_m^T.
\end{equation}
with
\begin{equation*}
\mat{S W}_{m + 1} = \mat{Q}_{m + 1} \mat{T}_{m + 1} = \begin{bmatrix} \mat{Q}_m &  \vec{q} \end{bmatrix} \begin{bmatrix} \mat{T}_m & \vec{t} \\ \vec{0} & t_{m + 1}  \end{bmatrix}
\end{equation*}
Note that $\mat{Q}_m^*  \mat{S} \vec{b} = \| \mat{S} \vec{b} \| \vec{e}_1$ is only valid for a positive $t_{11}$, which is not guaranteed by \texttt{dgeqrf}. Therefore, we compute $\text{sign}(t_{11}) \| \mat{S} \vec{b} \| \vec{e}_1$ instead of just $\| \mat{S} \vec{b} \| \vec{e}_1$ as all matrices are real in our experiments.

\textbf{Spectral deflation:} At each restart cycle, the algorithm first computes the Schur decomposition via \texttt{dgees}, reorders the Schur form so that the $\ell$ smallest eigenvalues appear at the start, and then constructs the deflation subspace with a single call to \texttt{gemm}, as shown in \cite{eiermann_deflated_2011}. The new basis vectors are then orthogonalized against the augmented subspace.

\end{appendix}

\section*{Competing interests}

All authors certify that they have no affiliations with or involvement in any organization or entity with any financial interest or non-financial interest in the subject matter or materials discussed in this manuscript.

\section*{Data availability}

The code for running the numerical examples can be found at \url{https://gitlab.com/nlg550/randomized-krylov}.

\bibliographystyle{spmpsci}
\bibliography{references}
\end{document}